\documentclass[12pt]{amsart}
\usepackage[margin=1in]{geometry}
\usepackage{bm}
\usepackage[dvipsnames,svgnames,x11names]{xcolor}
\usepackage[hyphens]{url}
\usepackage[colorlinks=true,
            linkcolor=Maroon,
            citecolor=blue,
            urlcolor=blue,
            hypertexnames=false,
            linktocpage]{hyperref}
\usepackage{bookmark}
\usepackage{amsmath,mathtools,amssymb,aliascnt}
\mathtoolsset{showonlyrefs=true}
\usepackage{graphicx}
\usepackage{tikz}
\usepackage{fancyhdr}
\usepackage{esint}
\usepackage{enumerate}
\usepackage{caption}
\allowdisplaybreaks

\newcommand{\tr}{\operatorname{tr}}

\newcommand{\Ric}{\operatorname{Ric}}

\theoremstyle{plain}
\newtheorem{theorem}{Theorem}[section]
\newtheorem*{theorem*}{Theorem}

\newaliascnt{proposition}{theorem}
\newtheorem{proposition}[proposition]{Proposition}
\aliascntresetthe{proposition}

\newaliascnt{lemma}{theorem}
\newtheorem{lemma}[lemma]{Lemma}
\aliascntresetthe{lemma}

\newaliascnt{corollary}{theorem}
\newtheorem{corollary}[corollary]{Corollary}
\aliascntresetthe{corollary}

\theoremstyle{definition}
\newaliascnt{definition}{theorem}
\newtheorem{definition}[definition]{Definition}
\aliascntresetthe{definition}

\theoremstyle{remark}
\newaliascnt{remark}{theorem}
\newtheorem{remark}[remark]{Remark}
\newtheorem*{remark*}{Remark}
\aliascntresetthe{remark}

\numberwithin{equation}{section}

\providecommand{\dd}{\,\mathrm{d}}
\providecommand{\diver}{\operatorname{div}}
\providecommand{\Hhyp}{\mathbb{H}}

\begin{document}

\title[A Reilly approach to weighted curvature integral inequalities]{A Reilly approach to weighted curvature integral inequalities in hyperbolic space}
\author[S. Bi]{Shengguang Bi}
\address{School of Mathematical Sciences, Beihang University, Beijing 100191, P.R. China}
\email{sy2509103@buaa.edu.cn}
\author[Y. Hu]{Yingxiang Hu}
\address{School of Mathematical Sciences, Beihang University, Beijing 100191, P.R. China}
\email{huyingxiang@buaa.edu.cn}
\thanks{This work was supported by the National Key Research and Development Program of China (Grant No.~2021YFA1001800) and the National Natural Science Foundation of China (Grant No.~12671060).}
\hypersetup{pdftitle={A Reilly approach to weighted curvature integral inequalities in hyperbolic space},pdfauthor={Shengguang Bi and Yingxiang Hu}}
\date{}

\begin{abstract}
We establish a polarized weighted Reilly identity and an exact deficit decomposition for weighted curvature integrals $I_j=\int_\Sigma VE_j\dd\mu$ in hyperbolic space, where $V=\cosh r$ and $E_j$ is the normalized $j$th mean curvature. For compact domains with smooth connected boundary, these yield $I_k^2\geq I_{k-1}I_{k+1}$ under intrinsic sub-staticity at $k=1$ and static convexity with $\kappa\in\Gamma_k\cap\overline{\Gamma}_{k+1}$ for $2\leq k\leq n-1$, with equality precisely for geodesic balls. Under $h\geq0$ and $Vh-V_\nu g\geq0$, log-concavity and a Euclidean mixed-volume representation recover the centered-ball comparisons conjectured by Hu and Li \cite[Conjecture~1.1]{HuLi}, allowing the weight center to lie outside the domain. Equality in these comparisons holds precisely for balls centered at the weight center.
\end{abstract}

\keywords{Reilly formula, weighted curvature integrals, static convexity, Brunn--Minkowski inequality, Alexandrov--Fenchel inequality}
\subjclass[2020]{53C42, 52A40}
\maketitle

\section{Introduction and main results}\label{sec:introduction}

Reilly's original formula \cite{Reilly} relates quadratic Hessian terms to Ricci curvature and boundary geometry. Ros \cite{Ros} applied this formula to the solution of a Dirichlet problem to prove the Heintze--Karcher inequality for compact Riemannian manifolds with nonnegative Ricci curvature and positive boundary mean curvature. Using higher-order Reilly formulas, Qiu--Xia \cite{QiuXiaNeumann} proved a family of Euclidean Alexandrov--Fenchel inequalities by solving classical Neumann problems for Hessian equations.

Auxiliary elliptic equations and the Bochner formula also underlie almost-Schur inequalities. De Lellis--Topping \cite{DeLellisTopping} proved a sharp bound for the scalar-curvature variance in terms of the trace-free Ricci tensor on closed connected manifolds of dimension at least three with nonnegative Ricci curvature. Perez \cite{Perez} obtained the corresponding mean-curvature estimate, involving the trace-free second fundamental form, for closed connected Euclidean hypersurfaces of dimension at least two with nonnegative Ricci curvature. For related extensions, higher-order estimates, and rigidity results, see \cite{GeWangXia,KwongAlmostSchur,ChengTensor,ChengZhou}.

Weighted extensions of Reilly's formula were established by Qiu--Xia \cite{QiuXiaWeighted} and Li--Xia \cite{LiXia}, with applications to Heintze--Karcher and Minkowski inequalities in space forms and sub-static manifolds. Li--Xia \cite{LiXia} combined the closed-manifold weighted Reilly formula with the Codazzi equation to bound the weighted variance of the mean curvature in terms of the trace-free second fundamental form. Their argument applies to closed connected hypersurfaces in hyperbolic space under intrinsic sub-staticity. Using the same identity, Zeng \cite{Zeng} obtained estimates for higher mean curvatures of closed connected sub-static hypersurfaces of dimension at least three in space forms. Huang, Ma, and Zhu \cite{HuangMaZhu} further extended the weighted Reilly formula to the drift Laplacian and obtained almost-Schur estimates for mean and higher mean curvatures. Related results for free boundary hypersurfaces in a ball were obtained in \cite{ChenHuLi}.

In this paper, we prove adjacent inequalities for weighted curvature integrals in hyperbolic space by extending the closed-hypersurface Reilly argument. We first establish a polarized weighted Reilly identity and then solve a linear elliptic equation involving the ratio of consecutive normalized mean curvatures. The auxiliary equation gives an exact deficit decomposition. The polarized identity and a pointwise mixed-discriminant inequality show that both terms are nonnegative under suitable pointwise tensor conditions. Their vanishing characterizes geodesic balls for the canonical static potential. A pointwise mixed-discriminant inequality and an integral Hessian identity are also used in the Bochner proof of the Alexandrov--Fenchel inequality by Shenfeld and van Handel \cite{ShenfeldVanHandel}.

Let $\Omega\subset\Hhyp^{n+1}$ be a compact domain with smooth connected boundary $\Sigma=\partial\Omega$, where $n\geq2$. Fix a point $p\in\Hhyp^{n+1}$, write $r=\operatorname{dist}(\cdot,p)$, and take $V=\cosh r$. We use $\bar g$ and $\bar\nabla$ to denote the hyperbolic metric and its Levi-Civita connection on $\Hhyp^{n+1}$. We also denote the induced metric by $g$, the outward unit normal by $\nu$, and the second fundamental form by $h$, with the convention $h(a,b)=\bar g(\bar\nabla_a\nu,b)$ for tangent vectors $a,b$. The principal curvatures $\kappa_1,\ldots,\kappa_n$ are the eigenvalues of the shape operator $g^{-1}h$. The (normalized) $j$th mean curvature is defined by
\begin{equation}\label{eq:notation}
E_j=\binom{n}{j}^{-1}\sigma_j(\kappa)=\binom{n}{j}^{-1}\sum_{1\leq i_1<\ldots<i_j\leq n} \kappa_{i_1}\cdots\kappa_{i_j},\quad 1\leq j\leq n,
\end{equation}
and $E_0=1$ by convention. We define the $j$th weighted curvature integral and the weighted volume by
\begin{equation}
I_j=\int_\Sigma VE_j\dd\mu,\qquad 0\leq j\leq n,
\end{equation}
and
\begin{equation}
I_{-1}=(n+1)\int_\Omega V\dd\mathrm{vol}. 
\end{equation}
Write $V_\nu=\bar g(\bar\nabla V,\nu)$. 
Using $g$ to identify symmetric two-tensors with self-adjoint endomorphisms, we denote the identity by $\mathrm{Id}$. Set $S=Vh-V_\nu\mathrm{Id}$, so $S\geq0$ is equivalent to $Vh-V_\nu g\geq0$. The Newton transformations of $h$ are
\begin{equation}\label{eq:newtontransformations} T_0=\mathrm{Id},\qquad T_j=\sigma_j\mathrm{Id}-hT_{j-1}\quad(1\leq j\leq n). \end{equation}
In a principal frame, $T_j$ has eigenvalues $\sigma_j(\kappa| i)$, where the $i$th argument is omitted, so in particular $T_n=0$.

The potential $V$ extends smoothly across $p$ and satisfies
\begin{equation}\label{eq:hyperbolicstaticpotential} 
\bar\nabla^2V=V\bar g,\qquad \bar\Delta V=(n+1)V,\qquad |\bar\nabla V|^2=V^2-1. 
\end{equation}
In particular, $|V_\nu|<V$. A domain is static convex if $S\geq0$ and strictly static convex if $S>0$.

For $1\leq\ell\leq n$, define the $\ell$th Gårding cone by
\begin{equation}\label{eq:Gardingcone} \Gamma_\ell=\{\kappa\in\mathbb R^n:\sigma_j(\kappa)>0\text{ for }1\leq j\leq\ell\}. \end{equation}
We write $\overline{\Gamma}_\ell$ for its closure.
The weighted Reilly argument yields the following adjacent inequalities.

\begin{theorem}\label{thm:allorders}
Let $\Omega\subset\Hhyp^{n+1}$ be a compact domain with smooth connected boundary, and let $V=\cosh\operatorname{dist}(\cdot,p)$. Fix an integer $1\leq k\leq n-1$ and assume the following.
\begin{enumerate}[(i)]
\item For $k=1$, assume only
\begin{equation}\label{eq:intrinsichypothesis} (\Delta V)g-\nabla^2V+V\Ric_\Sigma\geq0. \end{equation}
\item For $2\leq k\leq n-1$, assume
\begin{equation}\label{eq:globalhypotheses} Vh-V_\nu g\geq0,\qquad \kappa\in\Gamma_k\cap\overline{\Gamma}_{k+1}. \end{equation}
\end{enumerate}
Then
\begin{equation}\label{eq:missing} \int_\Sigma VE_{k-1}\left(\frac{E_k}{E_{k-1}}-\frac{I_k}{I_{k-1}}\right)^2\dd\mu\leq\int_\Sigma V\left(\frac{E_k^2}{E_{k-1}}-E_{k+1}\right)\dd\mu. \end{equation}
Equivalently,
\begin{equation}\label{eq:allorders} I_k^2\geq I_{k-1}I_{k+1}. \end{equation}
Equality holds if and only if $\Omega$ is a geodesic ball, whose center may differ from $p$.
\end{theorem}

The case $k=1$ recovers the weighted almost-Schur inequality under intrinsic sub-staticity. For the case $2\leq k\leq n-1$, inequality \eqref{eq:missing} estimates the variance of $E_k/E_{k-1}$ with respect to the measure $VE_{k-1}\dd\mu$. For ambient static potentials, \eqref{eq:Qfactor} identifies the intrinsic sub-static tensor with $ST_1$, and \eqref{eq:mixedReilly} shows that $ST_k$ is its natural higher-order counterpart in the polarized Reilly formula. The cone hypothesis guarantees ellipticity and, together with static convexity, supplies the required sign conditions while allowing $E_{k+1}=0$. In the case $k=0$, static convexity gives Xia's weighted-volume inequality \cite{Xia}. Outer parallel approximation then yields log-concavity of the positive sequence $I_{-1},I_0,\ldots,I_n$ under $h\geq0$ and $S\geq0$; see \autoref{cor:mixedpolynomial}.

A sharp lower bound for $I_1$ in terms of the unweighted area was proved by de Lima and Gir\~ao \cite{deLimaGirao} for star-shaped, strictly mean-convex hypersurfaces. The corresponding area inequalities for higher odd-order integrals were established under horospherical convexity in \cite{GeWangWu}. Within this class, Hu, Li, and Wei \cite{HuLiWei} used locally constrained curvature flows to compare weighted curvature integrals with unweighted quermassintegrals, extending the area inequalities to a larger family of sharp comparisons.

This leads to the problem of comparing weighted curvature integrals of different orders. Hu and Li \cite{HuLi} conjectured that, among static convex domains containing the weight center, centered balls minimize a higher-order weighted curvature integral when a lower-order one is prescribed. The static convexity condition in that work is $h\geq(V_\nu/V)g>0$. Combining the adjacent inequalities with a highest-order estimate from the mixed-volume representation, we recover these comparisons under the weaker assumptions below.

\begin{theorem}\label{thm:HuLiconjecture}
Let $\Omega\subset\Hhyp^{n+1}$ be a compact domain with smooth connected boundary, let $V=\cosh\operatorname{dist}(\cdot,p)$, and assume
\begin{equation}\label{eq:HuLihypotheses} h\geq0,\qquad Vh-V_\nu g\geq0. \end{equation}
The weight center $p$ is not required to lie in $\Omega$.
For integers $-1\leq k<\ell\leq n$, choose $R>0$ by $I_k(\Omega)=I_k(B_R(p))$. Then
\begin{equation}\label{eq:HuLiconjecture} I_\ell(\Omega)\geq I_\ell(B_R(p)). \end{equation}
Equality holds if and only if $\Omega=B_R(p)$.
\end{theorem}

To provide a geometric interpretation of the adjacent inequalities, we use spatial projection in the hyperboloid model centered at $p$. Writing $X=(V,x)$ and $\pi(X)=x$, we have
\begin{equation}\label{eq:canonicalprojection} \pi(\exp_p(r\theta))=\sinh r\,\theta,\qquad \theta\in\mathbb S^n. \end{equation}
Static convexity is equivalent to Euclidean convexity of the spatial projection \cite{HuLiWanXu}. The flat connection used below is a special case of the affine connections introduced in \cite{LiXiaAffine}, which give a Bochner interpretation of the weighted Reilly formula. If $h>0$ and $S>0$, we express the weighted curvature integrals as Euclidean mixed volumes of the spatial body and its projected normal body. These formulas extend by approximation to $h\geq0$ and $S\geq0$ and provide a generalized Brunn--Minkowski interpretation of the adjacent inequalities. Jensen's inequality and the Euclidean isoperimetric inequality yield the highest-order estimate, with equality precisely for balls centered at $p$. Brunn--Minkowski concavity also yields boundary Poincar\'e inequalities \cite{Colesanti}, which admit extensions to weighted Riemannian manifolds under curvature-dimension conditions via a generalized Reilly formula \cite{KolesnikovMilman}. The Euclidean case underlies the support variations in \autoref{sec:variations}.

In \autoref{sec:adjacent}, we prove the polarized weighted Reilly identity and the exact deficit decomposition, and establish the adjacent inequalities with rigidity. \autoref{sec:staticdomains} develops the projection geometry and the Brunn--Minkowski formulation. \autoref{sec:projection} gives the mixed-volume representation and proves \autoref{thm:HuLiconjecture}. Finally, \autoref{sec:variations} relates the deficit $I_0^2-I_{-1}I_1$ in the weighted-volume inequality to support variations and parallel domains.

\medskip
\noindent\textbf{Added note:} We notice that Xu \cite{XuAF} has recently established weighted Alexandrov--Fenchel inequalities in space forms by a mixed-volume representation. In the hyperbolic case, his results with $\varepsilon=0$ include the representation \eqref{eq:originalmixedvolumes}, the log-concavity of the positive sequence $I_{-1},I_0,\ldots,I_n$ under $h\geq0$ and $S\geq0$, and the comparisons in \autoref{thm:HuLiconjecture}. For the centered-ball comparisons, both proofs use a highest-order estimate from Jensen's inequality and the Euclidean isoperimetric inequality, followed by induction using log-concavity. Our proof of the adjacent inequalities proceeds through a linear auxiliary equation on $\Sigma$ and the polarized Reilly identity \eqref{eq:mixedReilly}, without invoking the mixed-volume representation. The exact deficit decomposition \eqref{eq:exactdeficit} gives the pointwise criterion and hyperbolic rigidity stated in \autoref{thm:pointwiseReilly} and \autoref{cor:hyperbolicrigidity}.

\section{The Reilly method and hyperbolic rigidity}\label{sec:adjacent}\label{sec:schur}

We extend the first-order argument in \cite[Sec.~6]{LiXia} using a polarized Reilly identity and an auxiliary elliptic equation involving the ratio $E_k/E_{k-1}$. We allow arbitrary positive boundary weights at the first order and positive ambient static potentials at higher orders. We use $S=Vh-V_\nu\mathrm{Id}$ and the Newton transformations \eqref{eq:newtontransformations}. The cones \eqref{eq:Gardingcone} are nested, with $\Gamma_n$ the positive orthant.

\subsection{Weighted volume, intrinsic sub-staticity, and variance}

For $V=\cosh\operatorname{dist}(\cdot,p)$, the weighted almost-Schur inequality \cite[Thm.~1.9]{LiXia} takes the form
\begin{equation}\label{eq:three} I_1^2\geq I_0I_2 \end{equation}
under horospherical convexity. In the case $k=0$, static convexity gives the weighted-volume inequality \cite[Thm.~1.1]{Xia}:
\begin{equation}\label{eq:volume} I_0^2\geq I_{-1}I_1. \end{equation}
Equality in either case holds precisely for geodesic balls, whose centers may differ from $p$.

The proof of \eqref{eq:volume} in \cite[Eqs.~(7), (10), and~(18)]{Xia} uses an auxiliary function satisfying
\begin{equation}\label{eq:XiaNeumann} \bar\Delta f-(n+1)f=1\ \text{in }\Omega,\qquad Vf_\nu-V_\nu f=\frac{I_{-1}}{(n+1)I_0}V\ \text{on }\Sigma. \end{equation}
The substitution $f=Vw$ gives a uniformly elliptic Neumann problem with compatible flux, hence \eqref{eq:XiaNeumann} has a smooth solution $f$, unique up to an additive multiple of $V$. With $Y=\nabla(f|_\Sigma)-(f|_\Sigma/V)\nabla V$, separating the trace-free part in the weighted Reilly formula gives
\begin{equation}\label{eq:XiaReillycomparison} \frac{nI_{-1}(I_0^2-I_{-1}I_1)}{(n+1)^2I_0^2}=\int_\Omega V\left|\bar\nabla^2f-\left(f+\frac{1}{n+1}\right)\bar g\right|^2\dd\mathrm{vol}+\int_\Sigma S(Y,Y)\dd\mu. \end{equation}
Both terms are nonnegative under static convexity. The argument below uses an auxiliary equation on $\Sigma$ to obtain the higher-order decomposition.

For the boundary formulas, we consider a smooth closed connected hypersurface $\Sigma$ in $\Hhyp^{n+1}$ with an arbitrary smooth positive weight $V$, and define
\begin{equation}\label{eq:Q} \mathcal Q=(\Delta V)g-\nabla^2V+V\Ric_\Sigma, \qquad \mathring{h}=h-E_1g. \end{equation}
The condition $\mathcal Q\geq0$ is intrinsic sub-staticity.
\begin{proposition}\label{prop:substaticfactor}
Suppose that $V$ extends to an ambient static potential satisfying $\bar\nabla^2V=V\bar g$. Identifying symmetric two-tensors with self-adjoint endomorphisms by $g$, we have
\begin{equation}\label{eq:Qfactor} \mathcal Q=ST_1=(Vh-V_\nu\mathrm{Id})T_1,\qquad T_1=(\tr h)\mathrm{Id}-h. \end{equation}
If $\kappa\in\Gamma_2$, then $T_1>0$, so intrinsic sub-staticity is equivalent to static convexity.
\end{proposition}

\begin{proof}
The restriction formulas are
\[ \nabla^2V=Vg-V_\nu h,\qquad \Delta V=nV-(\tr h)V_\nu. \]
Substitution of the Gauss equation $\Ric_\Sigma=-(n-1)g+(\tr h)h-h^2$ into \eqref{eq:Q} gives \eqref{eq:Qfactor}; see also \cite[Eq.~(6.2)]{LiXia}. If $\kappa\in\Gamma_2$, then $\tr h>0$ and $(\tr h)^2>|h|^2$, hence $\tr h>|h|\geq\kappa_i$ for every $i$. Thus $T_1>0$, and the equivalence follows because $S$ and $T_1$ commute.
\end{proof}

For the canonical potential, ambient staticity therefore does not by itself imply intrinsic sub-staticity: the boundary tensor is $ST_1$. Under horospherical convexity, both factors are positive definite.

The closed-manifold weighted Reilly formula \cite[Thm.~1.1]{LiXia} reads
\begin{equation}\label{eq:firstReilly} \int_\Sigma V\bigl[(\tr\mathcal A)^2-|\mathcal A|^2\bigr]\dd\mu=\int_\Sigma\mathcal Q(Y,Y)\dd\mu, \end{equation}
where $\mathcal A=\nabla^2f-(f/V)\nabla^2V$ and $Y=\nabla f-(f/V)\nabla V$. It holds for every $f\in C^\infty(\Sigma)$ and every smooth positive weight $V$ on $\Sigma$. Hence $\mathcal Q\geq0$ implies nonnegativity of the left side for every $f$.

The Codazzi and Reilly argument in \cite[Sec.~6]{LiXia} applies whenever $V>0$ and $\mathcal Q\geq0$, and yields
\begin{equation}\label{eq:schur} \int_\Sigma V\left(E_1-\frac{I_1}{I_0}\right)^2\dd\mu\leq\frac{1}{n(n-1)}\int_\Sigma V|\mathring{h}|^2\dd\mu. \end{equation}
Since
\begin{equation}\label{eq:newton2} E_1^2-E_2=\frac{|\mathring{h}|^2}{n(n-1)}, \end{equation}
the variance estimate \eqref{eq:schur} is equivalent to the integral inequality \eqref{eq:three}. We recover both at $k=1$ in \autoref{thm:pointwiseReilly}.

Whenever $V>0$ and $E_{k-1}>0$, set
\begin{equation}\label{eq:rho} e=\frac{E_k}{E_{k-1}}-\frac{I_k}{I_{k-1}}. \end{equation}
Expanding its weighted square yields
\begin{equation}\label{eq:generaldeficit} \frac{I_k^2-I_{k-1}I_{k+1}}{I_{k-1}}=\int_\Sigma V\left(\frac{E_k^2}{E_{k-1}}-E_{k+1}\right)\dd\mu-\int_\Sigma VE_{k-1}e^2\dd\mu. \end{equation}

By \eqref{eq:generaldeficit}, the adjacent integral inequality \eqref{eq:allorders} is equivalent to the variance estimate \eqref{eq:missing}. If $\kappa\in\Gamma_k$, the pointwise quantity $E_k^2/E_{k-1}-E_{k+1}$ is nonnegative: for $E_{k+1}>0$, this follows from Newton's inequality, and for $E_{k+1}\leq0$, it is immediate. At $k=1$, it follows from \eqref{eq:newton2}. It remains to estimate the weighted variance. For $k\geq2$, the hypotheses \eqref{eq:globalhypotheses} ensure $ST_k\geq0$ and the bilinear-form sign condition in \eqref{eq:directReillyconditions}, as verified in the proof of \autoref{thm:allorders}.

\subsection{The bilinear form and polarized Reilly identity}\label{subsec:polarizedReilly}

Let $D$ be the mixed discriminant normalized by $D(A[n])=\det A$, where brackets indicate repeated arguments. For $1\leq k\leq n-1$, define
\begin{equation}\label{eq:polarizedform} \mathcal B_k(A,B)=D(h[k-1],A,B,\mathrm{Id}[n-k-1]). \end{equation}
Intrinsic indices are raised and lowered with $g$ and repeated indices are summed from $1$ to $n$. Components of a self-adjoint endomorphism refer to its associated covariant tensor, so $A^{ij}=g^{ip}g^{jq}A_{pq}$.

For the mixed discriminant and the form \eqref{eq:polarizedform}, we have $E_j(h)=D(h[j],\mathrm{Id}[n-j])$ and
\begin{equation}\label{eq:polarizedvalues} \mathcal B_k(\mathrm{Id},\mathrm{Id})=E_{k-1},\qquad \mathcal B_k(h,\mathrm{Id})=E_k,\qquad \mathcal B_k(h,h)=E_{k+1}. \end{equation}

Set
\begin{equation}\label{eq:PR} P_k=\frac{T_{k-1}}{k\binom{n}{k}}. \end{equation}
Then $\tr P_k=E_{k-1}$ and
\begin{equation}\label{eq:polarizedlinear} \mathcal B_k(\mathrm{Id},A)=P_k^{ij}A_{ij},\qquad \mathcal B_k(h,A)=\frac{T_k^{ij}A_{ij}}{(k+1)\binom{n}{k+1}}. \end{equation}
The Codazzi equations in $\Hhyp^{n+1}$ imply $\diver T_j=0$ and hence $\diver P_k=0$.

\begin{theorem}[Polarized weighted Reilly identity]\label{lem:mixedReilly}
Let $\Sigma$ be a smooth closed hypersurface in $\Hhyp^{n+1}$, and let $V$ be a smooth function in a neighborhood of $\Sigma$ with $V>0$ on $\Sigma$ and $\bar\nabla^2V=V\bar g$. For $1\leq k\leq n-1$ and $f\in C^\infty(\Sigma)$, set
\begin{equation}\label{eq:AYhigher} \mathcal A=\nabla^2f-\frac fV\nabla^2V,\qquad Y=\nabla f-\frac fV\nabla V. \end{equation}
Then
\begin{equation}\label{eq:mixedReilly} \int_\Sigma V\mathcal B_k(\mathcal A,\mathcal A)\dd\mu =\frac{1}{(k+1)\binom{n}{k+1}}\int_\Sigma(ST_k)(Y,Y)\dd\mu. \end{equation}
Since $S$ commutes with $T_k$, the condition
\begin{equation}\label{eq:Reillysign} ST_k\geq0 \end{equation}
ensures that the integral on the left of \eqref{eq:mixedReilly} is nonnegative for every $f$. At $k=1$, \eqref{eq:Reillysign} is exactly intrinsic sub-staticity because $ST_1=\mathcal Q$ by \eqref{eq:Qfactor}. Thus $ST_k\geq0$ is a natural higher-order analogue of this sign condition in the polarized identity \eqref{eq:mixedReilly}. In a principal frame, this condition is $(V\kappa_i-V_\nu)\sigma_k(\kappa| i)\geq0$ for every $i$. When $T_k>0$, it is equivalent to static convexity. Under $\kappa\in\Gamma_k\cap\overline{\Gamma}_{k+1}$, we have $T_k\geq0$ by continuity, so static convexity is sufficient even when $T_k$ is singular.
\end{theorem}

\begin{proof}
Work in a local orthonormal frame and sum over repeated indices. Write $V_i=\nabla_iV$ and $V_{ij}=\nabla_i\nabla_jV$. Define
\begin{equation}\label{eq:Ctensor} \mathcal C_k^{ij\ell q} =\frac{ \delta^{i_1\cdots i_{k-1}i\ell}_{j_1\cdots j_{k-1}jq} h_{i_1}^{j_1}\cdots h_{i_{k-1}}^{j_{k-1}} }{(k+1)!\binom{n}{k+1}}, \end{equation}
where $\delta$ is the generalized Kronecker delta. For $k=1$, the product is understood to be $1$. The determinant polarization is
\[ \mathcal B_k(A,B)=\mathcal C_k^{ij\ell q}A_{ij}B_{\ell q}. \]
The tensor $\mathcal C_k$ is antisymmetric in $i,\ell$ and in $j,q$. By the Codazzi equation, $\nabla_i h_{i_s}{}^{j_s}=\nabla_{i_s}h_i{}^{j_s}$ for each factor. Its contraction with the generalized delta vanishes by antisymmetry in $i,i_s$. Hence $\nabla_i\mathcal C_k^{ij\ell q}=0$. The identity $\mathcal C_k^{jiq\ell}=\mathcal C_k^{ij\ell q}$ also shows that $\mathcal C_k^{ij\ell q}\mathcal A_{\ell q}$ is symmetric in $i,j$.

Using $\mathcal A_{ij}=\nabla_iY_j+Y_iV_j/V$, the symmetry above reduces the integrand to
\[ V\mathcal C_k^{ij\ell q}\mathcal A_{\ell q}\mathcal A_{ij} =\mathcal C_k^{ij\ell q}\mathcal A_{\ell q}\nabla_i(VY_j). \]
Integrating by parts, we obtain
\begin{equation}\label{eq:mixedIBP} \int_\Sigma V\mathcal B_k(\mathcal A,\mathcal A)\dd\mu =-\int_\Sigma VY_j\mathcal C_k^{ij\ell q}\nabla_i\mathcal A_{\ell q}\dd\mu. \end{equation}
Our convention is $R_{i\ell q}{}^pe_p=R(e_i,e_\ell)e_q$, with $R_{ijj}{}^i=\kappa_i\kappa_j-1$ for $i\ne j$ in a principal frame. Commuting covariant derivatives of $f$ and $V$ and using $\nabla_i(f/V)=Y_i/V$, we obtain
\begin{equation}\label{eq:Acommutator} \nabla_i\mathcal A_{\ell q}-\nabla_\ell\mathcal A_{iq} =-R_{i\ell q}{}^pY_p-\frac{Y_i}{V}V_{\ell q}+\frac{Y_\ell}{V}V_{iq}. \end{equation}
Using the antisymmetry in $i,\ell$, we rewrite \eqref{eq:mixedIBP} as
\begin{equation}\label{eq:mixedcurvature} \int_\Sigma V\mathcal B_k(\mathcal A,\mathcal A)\dd\mu =\int_\Sigma\left[\frac V2\mathcal C_k^{ij\ell q}R_{i\ell q}{}^pY_jY_p+\mathcal C_k^{ij\ell q}V_{\ell q}Y_iY_j\right]\dd\mu. \end{equation}

At a point where $h$ is diagonal, the ambient static equation and the Gauss equations read
\[ V_{jj}=V-V_\nu\kappa_j,\qquad R_{ijj}{}^i=\kappa_i\kappa_j-1\quad(i\ne j). \]
For $i\ne j$, set
\[ c_{ij}=\mathcal C_k^{ii jj}=\frac{\sigma_{k-1}(\kappa| i,j)}{k(k+1)\binom{n}{k+1}},\qquad \mathcal C_k^{ijji}=-c_{ij},\qquad c_{ij}=c_{ji}, \]
where $\kappa| i,j$ denotes the principal curvatures with indices $i,j$ omitted. The two contractions in \eqref{eq:mixedcurvature} are
\[
\begin{aligned}
\frac12\mathcal C_k^{ij\ell q}R_{i\ell q}{}^pY_jY_p&=\sum_{i\ne j}c_{ij}(\kappa_i\kappa_j-1)Y_i^2,\\
\mathcal C_k^{ij\ell q}V_{\ell q}Y_iY_j&=\sum_{i\ne j}c_{ij}(V-V_\nu\kappa_j)Y_i^2.
\end{aligned}
\]
Since $V(\kappa_i\kappa_j-1)+V-V_\nu\kappa_j=(V\kappa_i-V_\nu)\kappa_j$, we have
\[ \int_\Sigma V\mathcal B_k(\mathcal A,\mathcal A)\dd\mu=\int_\Sigma\sum_{i\ne j}c_{ij}(V\kappa_i-V_\nu)\kappa_jY_i^2\dd\mu. \]
Finally, the identity $\sum_{j\ne i}\kappa_j\sigma_{k-1}(\kappa| i,j)=k\sigma_k(\kappa| i)$ implies
\[ \sum_{j\ne i}c_{ij}\kappa_j=\frac{(T_k)_{ii}}{(k+1)\binom{n}{k+1}}. \]
Substitution proves \eqref{eq:mixedReilly}.
\end{proof}

\begin{remark}
For $k=1$, the left side of \eqref{eq:mixedReilly} is
\[ \frac{1}{n(n-1)}\int_\Sigma V\bigl[(\tr\mathcal A)^2-|\mathcal A|^2\bigr]\dd\mu, \]
and the right side is $\int_\Sigma\mathcal Q(Y,Y)\dd\mu/[n(n-1)]$. Thus the identity reduces to the closed-manifold case of \cite[Thm.~1.1]{LiXia}.

For $V=\cosh r$ on a hyperbolic geodesic sphere centered at $p$, the potential and the principal curvature $\kappa$ are constant, so the Bochner identity and $V\kappa-V_\nu=V(\kappa^2-1)/\kappa$ reduce both sides to
\[ \frac{V\kappa^{k-1}(\kappa^2-1)}{n} \int_\Sigma|\nabla f|^2\dd\mu. \]
\end{remark}

\begin{lemma}\label{lem:matrix}
For $k=1$, let $h$ be any symmetric endomorphism. For $2\leq k\leq n-1$, assume $\kappa(h)\in\Gamma_k\cap\overline{\Gamma}_{k+1}$. On the subspace $\{Z:\mathcal B_k(\mathrm{Id},Z)=0\}$, the form $-\mathcal B_k$ is positive semidefinite.
If $h>0$, the restricted form is positive definite.
\end{lemma}

\begin{proof}
For $k=1$,
\[ \mathcal B_1(A,B)=\frac{(\tr A)(\tr B)-\tr(AB)}{n(n-1)}. \]
Thus $\mathcal B_1(\mathrm{Id},Z)=\tr Z/n$, and on its kernel $-\mathcal B_1(Z,Z)=|Z|^2/[n(n-1)]\geq0$.

For $k\geq2$, first assume $E_{k+1}(h)>0$, so $\kappa(h)\in\Gamma_{k+1}$. On the space of real symmetric matrices, $E_{k+1}^{1/(k+1)}$ is concave on the open matrix cone $\{A\in\operatorname{Sym}(n):\kappa(A)\in\Gamma_{k+1}\}$ by G\aa rding's theorem \cite[Cors.~2.13(d) and~2.23]{HarveyLawson}. Thus $h+tZ$ remains in this cone for sufficiently small $|t|$, for every symmetric matrix $Z$. The first two derivatives are $\mathrm dE_{k+1}(h)[Z]=(k+1)\mathcal B_k(h,Z)$ and $\mathrm d^2E_{k+1}(h)[Z,Z]=k(k+1)\mathcal B_k(Z,Z)$. Differentiating the concavity inequality at $h$, we obtain
\begin{equation}\label{eq:GardingHessian} E_{k+1}\mathcal B_k(Z,Z)\leq\mathcal B_k(h,Z)^2. \end{equation}
Since $E_{k+1}>0$, the form $\mathcal B_k$ is nonpositive on $\{Z:\mathcal B_k(h,Z)=0\}$ and has at most one positive eigenvalue. On the other hand, $\mathcal B_k(\mathrm{Id},\mathrm{Id})=E_{k-1}>0$. If $Z$ were orthogonal to $\mathrm{Id}$ for $\mathcal B_k$ and satisfied $\mathcal B_k(Z,Z)>0$, the restriction of $\mathcal B_k$ to $\operatorname{span}\{\mathrm{Id},Z\}$ would be positive definite, a contradiction. Thus $-\mathcal B_k$ is positive semidefinite on the stated subspace.

If $E_{k+1}(h)=0$, use the identity
\begin{equation}\label{eq:coneshift} \sigma_j(\kappa+\eta\mathbf1)=\sum_{\ell=0}^j\binom{n-\ell}{j-\ell}\eta^{j-\ell}\sigma_\ell(\kappa),\qquad \sigma_0=1,\quad\eta>0. \end{equation}
For $j\leq k+1$, every term is nonnegative and at least one is positive, because $\kappa\in\Gamma_k$ and $\sigma_{k+1}=0$. Hence $h+\eta\mathrm{Id}\in\Gamma_{k+1}$ for every $\eta>0$. By continuity, the limiting bilinear form has at most one positive eigenvalue. Since $\mathcal B_k(\mathrm{Id},\mathrm{Id})=E_{k-1}>0$, the same two-dimensional argument proves nonpositivity on the stated kernel.

When $h>0$, all the fixed arguments $h[k-1],\mathrm{Id}[n-k-1]$ are positive definite. If $\mathcal B_k(\mathrm{Id},Z)=0$ and $\mathcal B_k(Z,Z)=0$, the mixed-discriminant inequality with arguments $Z,\mathrm{Id}$ is an equality. The equality characterization for positive definite fixed arguments applies to an arbitrary symmetric $Z$ and implies that $Z$ is a scalar multiple of $\mathrm{Id}$ \cite[Thm.~1.3 and Sec.~6.2]{ShenfeldVanHandel}. Since $\mathcal B_k(\mathrm{Id},\mathrm{Id})=E_{k-1}>0$, the orthogonality condition makes that scalar zero. Thus $-\mathcal B_k$ is positive definite on the indicated kernel.
\end{proof}

\subsection{The exact deficit identity and pointwise conditions}\label{subsec:adjacentproofs}

\begin{proposition}[Exact deficit decomposition]\label{prop:exactdeficit}
Let $\Sigma$ be smooth, closed, and connected in $\Hhyp^{n+1}$, and let $V\in C^\infty(\Sigma)$ be positive. Fix $1\leq k\leq n-1$ and assume $P_k>0$. Put $c=I_k/I_{k-1}$ and let $f$ be the unique smooth solution of
\begin{equation}\label{eq:newtonPDE} P_k^{ij}\left(\nabla_i\nabla_jf-\frac fV\nabla_i\nabla_jV\right)=E_{k-1}e,\qquad \int_\Sigma fV\dd\mu=0. \end{equation}
With $\mathcal A,Y$ as in \eqref{eq:AYhigher}, set $W=h-\mathcal A-c\mathrm{Id}$ and $J_k=\int_\Sigma V\mathcal B_k(\mathcal A,\mathcal A)\dd\mu$. Then $\mathcal B_k(\mathrm{Id},W)=0$ and
\begin{equation}\label{eq:exactdeficit} \frac{I_k^2-I_{k-1}I_{k+1}}{I_{k-1}}=J_k-\int_\Sigma V\mathcal B_k(W,W)\dd\mu. \end{equation}
At $k=1$, for any positive weight,
\begin{equation}\label{eq:exactfirstReilly} J_1=\frac1{n(n-1)}\int_\Sigma\mathcal Q(Y,Y)\dd\mu. \end{equation}
If $V$ is an ambient static potential, then at every order
\begin{equation}\label{eq:exactstaticReilly} J_k=\frac1{(k+1)\binom n{k+1}}\int_\Sigma(ST_k)(Y,Y)\dd\mu. \end{equation}
\end{proposition}

\begin{proof}
Since $\tr P_k=E_{k-1}>0$, the ratio and its weighted mean are well defined. The Codazzi equations imply $\diver P_k=0$, and $f=Vw$ transforms \eqref{eq:newtonPDE} into
\begin{equation}\label{eq:newtonPDEdiv} \diver(V^2P_k\nabla w)=VE_{k-1}e. \end{equation}
Since $\Sigma$ is compact, $V>0$ and $P_k>0$ imply $V^2P_k\geq a_0\mathrm{Id}$ for some $a_0>0$, so the transformed equation is uniformly elliptic. By the definition of $e$, its compatibility condition is
\[
\int_\Sigma VE_{k-1}e\dd\mu=0.
\]
The transformed operator is self-adjoint in $L^2(\Sigma,\dd\mu)$. Integration by parts for a homogeneous solution yields $\int_\Sigma V^2P_k(\nabla w,\nabla w)\dd\mu=0$, so connectedness makes $w$ constant. The Fredholm alternative and elliptic regularity therefore ensure a smooth solution, unique up to an additive constant. Equivalently, the original self-adjoint operator has kernel $\operatorname{span}\{V\}$, and the compatibility condition above is orthogonality to this kernel. The normalization $\int_\Sigma V^2w\dd\mu=0$ fixes the constant uniquely because $\int_\Sigma V^2\dd\mu>0$. The quantities $\mathcal A,Y,W$ and $J_k$ are invariant under $f\mapsto f+\alpha V$.

By \eqref{eq:polarizedlinear} and the auxiliary equation, $\mathcal B_k(\mathrm{Id},\mathcal A)=E_{k-1}e$, and hence $\mathcal B_k(\mathrm{Id},W)=E_k-E_{k-1}(e+c)=0$. The divergence-free Newton transformation also gives
\begin{equation}\label{eq:newtondivergence} \diver\bigl(T_k(V\nabla f-f\nabla V)\bigr)=V T_k^{ij}\mathcal A_{ij}. \end{equation}
Integration, \eqref{eq:polarizedlinear}, and the compatibility condition imply
\[ \int_\Sigma V\mathcal B_k(h,\mathcal A)\dd\mu=0,\qquad \int_\Sigma V\mathcal B_k(\mathrm{Id},\mathcal A)\dd\mu=0. \]
Expanding $W=h-c\mathrm{Id}-\mathcal A$, we obtain
\[ \int_\Sigma V\mathcal B_k(W,W)\dd\mu=I_{k+1}-2cI_k+c^2I_{k-1}+J_k=J_k-\frac{I_k^2-I_{k-1}I_{k+1}}{I_{k-1}}. \]
This proves \eqref{eq:exactdeficit}. Equation~\eqref{eq:exactfirstReilly} follows from \eqref{eq:firstReilly}, and \eqref{eq:exactstaticReilly} follows from \autoref{lem:mixedReilly}.
\end{proof}

\begin{theorem}\label{thm:pointwiseReilly}
Let $\Omega\subset\Hhyp^{n+1}$ be a compact domain with smooth connected boundary $\Sigma$, and let $V\in C^\infty(\Sigma)$ be positive. Fix $1\leq k\leq n-1$ and assume the following.
\begin{enumerate}[(i)]
\item For $k=1$, assume only \eqref{eq:intrinsichypothesis}.
\item For $2\leq k\leq n-1$, assume that $V$ extends smoothly to a neighborhood of $\Sigma$ and satisfies $\bar\nabla^2V=V\bar g$. At every point, require
\begin{equation}\label{eq:directReillyconditions} T_{k-1}>0,\qquad ST_k\geq0,\qquad \mathcal B_k(Z,Z)\leq0. \end{equation}
The last inequality is required for every symmetric $Z$ with $(T_{k-1})^{ij}Z_{ij}=0$.
\end{enumerate}
Then the variance estimate \eqref{eq:missing}, equivalently \eqref{eq:allorders}, holds. If $k=1$ and the tensor in \eqref{eq:intrinsichypothesis} is positive definite, equality holds if and only if $\Sigma$ is a geodesic sphere. For $V=\cosh\operatorname{dist}(\cdot,p)$, \autoref{cor:hyperbolicrigidity} characterizes equality under the stated hypotheses at every order.
\end{theorem}

\begin{proof}[Proof of \autoref{thm:pointwiseReilly}]
For $k=1$, $P_1=\mathrm{Id}/n>0$. For $k\geq2$, the first condition in \eqref{eq:directReillyconditions} implies $P_k>0$. \autoref{prop:exactdeficit} therefore applies. Its tensor $W$ satisfies $\mathcal B_k(\mathrm{Id},W)=0$. At $k=1$, $-\mathcal B_1(W,W)=|W|^2/[n(n-1)]$ by the explicit formula for $\mathcal B_1$, and $J_1\geq0$ by $\mathcal Q\geq0$. At higher orders, the other two conditions imply $J_k\geq0$ and $-\mathcal B_k(W,W)\geq0$. The two terms in \eqref{eq:exactdeficit} are nonnegative, proving \eqref{eq:allorders}. Equation~\eqref{eq:generaldeficit} proves the equivalent variance estimate.

If $k=1$, $\mathcal Q>0$, and equality holds, then \eqref{eq:exactdeficit} and \eqref{eq:exactfirstReilly} imply $Y=0$ and $W=0$. Since $Y=V\nabla(f/V)$, we have $\mathcal A=0$ and hence $h=c\mathrm{Id}$ for the constant $c=I_1/I_0$. Comparison with an enclosing geodesic ball implies $c>1$, so the connected boundary is a geodesic sphere. Conversely, its constant principal curvature implies $I_j=c^jI_0$ for $0\leq j\leq n$, and equality follows.
\end{proof}

\subsection{Rigidity and interpolation}\label{subsec:rigidity}

For the canonical potential $V=\cosh\operatorname{dist}(\cdot,p)$, let $R>0$, $q\in\Hhyp^{n+1}$, and $d=\operatorname{dist}(p,q)$. Then
\begin{equation}\label{eq:offcenterballcalibration} I_j(B_R(q))=\cosh d\,\omega_n\sinh^{n-j} R\cosh^{j+1} R,\qquad -1\leq j\leq n. \end{equation}
Here $\omega_n=|\mathbb S^n|$ is the area of the unit $n$-sphere.
Indeed, choose the polar axis about $q$ towards $p$ when $d>0$, and arbitrarily when $d=0$. The hyperbolic cosine rule reads $V(\exp_q(s\theta))=\cosh d\cosh s-\sinh d\sinh s\,\theta_1$. On $\partial B_R(q)$, the outward normal is $\partial_s$, so
\[
V_\nu=\cosh d\sinh R-\sinh d\cosh R\,\theta_1,\qquad V\coth R-V_\nu=\frac{\cosh d}{\sinh R}.
\]
Thus $h=(\coth R)g>g$ and $S=(\cosh d/\sinh R)\mathrm{Id}>0$ for every ball. The term containing $\theta_1$ has zero spherical integral. For $j\geq0$, use $E_j=\coth^j R$ and $\dd\mu=\sinh^n R\dd\mu_{\mathbb S}$. For $j=-1$, integrate $(n+1)\omega_n\cosh d\cosh s\,\sinh^n s$ from $0$ to $R$. This proves \eqref{eq:offcenterballcalibration}, including the case $j=-1$.

\begin{lemma}\label{lem:hyperboliccontinuation}
Let $\Omega\subset\Hhyp^{n+1}$ be a compact domain with smooth connected boundary, and let $V=\cosh\operatorname{dist}(\cdot,p)$ and $S=Vh-V_\nu\mathrm{Id}$. The set
\[ \mathcal U=\{x\in\Sigma:h(x)>0,\ S(x)>0\} \]
is nonempty and open. If $h=c\mathrm{Id}$ on $\mathcal U$ for a single constant $c$, then $c>1$ and $\Omega$ is a geodesic ball.
\end{lemma}

\begin{proof}
Let $R_*=\max_{\overline\Omega}r>0$. A maximum point $x_*$ lies on $\Sigma$, since $r$ has no critical point away from $p$. The inclusion $\Omega\subset\overline B_{R_*}(p)$ implies $\nu=\partial_r$ at $x_*$. There $0\geq\nabla^2r=\coth R_*\,g-h$, and hence
\begin{equation}\label{eq:canonicalellipticpoint} h(x_*)\geq\coth R_*\,\mathrm{Id}>\mathrm{Id},\qquad S(x_*)\geq\frac1{\sinh R_*}\mathrm{Id}>0. \end{equation}
The second inequality uses $V=\cosh R_*$ and $V_\nu=\sinh R_*$ at $x_*$. Thus $x_*\in\mathcal U$. If $h=c\mathrm{Id}$ on $\mathcal U$ for the same constant on every component, then $c\geq\coth R_*>1$. Let $x_j\in\mathcal U$ converge to $x_\infty\in\Sigma$. In a smooth local frame, $h(x_\infty)=c\mathrm{Id}$ by continuity. The canonical potential satisfies $|V_\nu|<V$, so
\[ S(x_\infty)=\bigl(V(x_\infty)c-V_\nu(x_\infty)\bigr)\mathrm{Id}\geq V(x_\infty)(c-1)\mathrm{Id}>0. \]
Thus $x_\infty\in\mathcal U$, which proves that $\mathcal U$ is closed as well as open, and connectedness of $\Sigma$ then implies $\mathcal U=\Sigma$.

In the hyperboloid model, write $\langle\cdot,\cdot\rangle_L$ for the Lorentzian inner product. The ambient Weingarten formula makes $\xi=\nu-cX$ constant on the connected boundary, with $\langle\xi,\xi\rangle_L=1-c^2<0$ and $\langle X,\xi\rangle_L=c$. Its time component is $V_\nu-cV<0$, so $q_0=-\xi/\sqrt{c^2-1}$ lies in $\Hhyp^{n+1}$ and
\[ -\langle X,q_0\rangle_L=\frac c{\sqrt{c^2-1}}=\cosh R_0,\qquad \coth R_0=c,\quad R_0>0. \]
Hence $\Sigma\subset\partial B_{R_0}(q_0)$. The inclusion has invertible differential between the tangent spaces of these $n$-dimensional hypersurfaces, so its image is open by the inverse function theorem and closed by compactness. The geodesic sphere is connected, so $\Sigma=\partial B_{R_0}(q_0)$. Since $\Omega$ is compact, it is the enclosed ball.
\end{proof}

The equality case of \eqref{eq:volume} under $S\geq0$ also follows from \autoref{lem:hyperboliccontinuation}. Suppose equality holds, let $f$ solve \eqref{eq:XiaNeumann}, and set $c_0=I_{-1}/((n+1)I_0)>0$. Equation~\eqref{eq:XiaReillycomparison} gives $\bar\nabla^2f=(f+1/(n+1))\bar g$ and $S(Y,Y)=0$. On $\mathcal U=\{h>0,S>0\}$, we have $Y=0$, so locally $f|_\Sigma=\alpha V$ for a constant $\alpha$. The boundary condition gives $(f-\alpha V)_\nu=c_0$. Restricting the ambient Hessian equation to this patch and using $\nabla^2(f|_\Sigma-\alpha V)=0$, we obtain
\[ 0=\frac{1}{n+1}g-c_0h. \]
Thus $h=(I_0/I_{-1})\mathrm{Id}$ throughout $\mathcal U$, with the same constant on every component. The lemma implies that $\Omega$ is a geodesic ball. Conversely, \eqref{eq:offcenterballcalibration} gives equality for every geodesic ball.

\begin{corollary}\label{cor:hyperbolicrigidity}\label{cor:canonicalrigidity}\label{cor:firstorderrigidity}
Let $\Omega\subset\Hhyp^{n+1}$ be a compact domain with smooth connected boundary, and let $V=\cosh\operatorname{dist}(\cdot,p)$. Fix $1\leq k\leq n-1$. At $k=1$, assume $\mathcal Q\geq0$. At $k\geq2$, assume \eqref{eq:directReillyconditions}, including $T_{k-1}>0$. Then equality in \eqref{eq:allorders}, equivalently in \eqref{eq:missing}, holds if and only if $\Omega$ is a geodesic ball, whose center may differ from $p$.
\end{corollary}

\begin{proof}
We adapt the continuation argument in \cite[Sec.~6]{LiXia} to the exact decomposition. The hypotheses ensure $P_k>0$, with $P_1=\mathrm{Id}/n$. Let $f$, $\mathcal A,Y,W$, and $c=I_k/I_{k-1}$ be as in \autoref{prop:exactdeficit}. Equality in the adjacent inequality \eqref{eq:allorders} makes the left side of \eqref{eq:exactdeficit} vanish. Since its two right-hand terms are nonnegative, $J_k=0$ and $\int_\Sigma V[-\mathcal B_k(W,W)]\dd\mu=0$. Continuity, $V>0$, and the pointwise sign condition imply $\mathcal B_k(W,W)=0$ everywhere. By \eqref{eq:exactfirstReilly} at $k=1$ and \eqref{eq:exactstaticReilly} at higher orders, $J_k=0$ also implies $\mathcal Q(Y,Y)=0$ or $(ST_k)(Y,Y)=0$, respectively, since these integrands are continuous and nonnegative.

Let $\mathcal U=\{h>0,\ S>0\}$, which is nonempty and open by \autoref{lem:hyperboliccontinuation}. On this set, $h>0$ implies $T_k>0$ for $1\leq k\leq n-1$. Since $S$ and $T_k$ commute, $ST_k>0$, with $\mathcal Q=ST_1>0$ at $k=1$. The vanishing of $(ST_k)(Y,Y)$ therefore forces $Y=0$ throughout $\mathcal U$. Thus $\nabla Y=0$ there, and $\mathcal A_{ij}=\nabla_iY_j+Y_iV_j/V$ implies $\mathcal A=0$. Since $\mathcal B_k(\mathrm{Id},W)=0$, the positive definiteness of $-\mathcal B_k$ on this subspace, given by \autoref{lem:matrix}, and $\mathcal B_k(W,W)=0$ imply $W=0$. Consequently $h=c\mathrm{Id}$ on all of $\mathcal U$, with the same global constant $c$ on every component. \autoref{lem:hyperboliccontinuation} proves that $\Omega$ is a geodesic ball.

Conversely, every geodesic ball has $h>0$ and $S>0$ by the calculation following \eqref{eq:offcenterballcalibration}, and all $T_j$ with $0\leq j\leq n-1$ are positive definite. \autoref{lem:matrix} verifies the bilinear-form condition, and \eqref{eq:offcenterballcalibration} implies equality.
\end{proof}

\begin{proof}[Proof of \autoref{thm:allorders}]
For $k=1$, apply \autoref{thm:pointwiseReilly}. For $k\geq2$, $\kappa\in\Gamma_k$ implies $T_{k-1}>0$. If $E_{k+1}>0$, then $\kappa\in\Gamma_{k+1}$ and $T_k>0$. At $E_{k+1}=0$, $T_k\geq0$ follows by approximation with $h+\eta\mathrm{Id}$. Since $S$ and $T_k$ commute, their nonnegativity implies $ST_k\geq0$. \autoref{lem:matrix} verifies the bilinear-form condition. The inequalities then follow from \autoref{thm:pointwiseReilly}, and the equality statement follows from \autoref{cor:hyperbolicrigidity}.
\end{proof}

\begin{remark}
For $n\geq3$, the inequality $I_2^2\geq I_1I_3$ holds under $\kappa\in\Gamma_2\cap\overline{\Gamma}_3$ and $Vh-V_\nu g\geq0$, with $V$ a positive ambient static potential. The equation and bilinear form are
\[ P_2=\frac{(\tr h)\mathrm{Id}-h}{n(n-1)},\qquad P_2^{ij}\mathcal A_{ij}=E_2-\frac{I_2}{I_1}E_1, \]
and
\begin{equation}\label{eq:B2explicit} \mathcal B_2(A,A) =\frac{(\tr h)\bigl[(\tr A)^2-|A|^2\bigr]-2(\tr A)\tr(hA)+2\tr(hA^2)}{n(n-1)(n-2)}. \end{equation}
The polarized Reilly identity takes the form
\[ \int_\Sigma V\mathcal B_2(\mathcal A,\mathcal A)\dd\mu =\frac{2}{n(n-1)(n-2)}\int_\Sigma\bigl((Vh-V_\nu\mathrm{Id})T_2\bigr)(Y,Y)\dd\mu\geq0. \]
Combining this nonnegativity with \eqref{eq:exactdeficit} and \autoref{lem:matrix}, we obtain $I_2^2\geq I_1I_3$. For $k=1$, the argument requires only the intrinsic condition $\mathcal Q\geq0$ and recovers \eqref{eq:schur}.
\end{remark}

\begin{corollary}\label{cor:geometry}
Let $\Omega\subset\Hhyp^{n+1}$ be a compact domain with smooth connected boundary, let $V=\cosh\operatorname{dist}(\cdot,p)$, and assume $\kappa_i\geq1$. Then \eqref{eq:allorders} holds at every $1\leq k\leq n-1$, with equality if and only if $\Omega$ is a geodesic ball.
\end{corollary}

\begin{proof}
Since $\kappa_i\geq1$ and $|V_\nu|<V$, we have $h>0$ and $S>0$. \autoref{thm:allorders} applies at every order.
\end{proof}

\begin{corollary}\label{cor:interpolation}
Let $\Omega\subset\Hhyp^{n+1}$ be a compact domain with smooth connected boundary, let $V=\cosh\operatorname{dist}(\cdot,p)$, and fix $1\leq\ell\leq n$. Assume
\begin{equation}\label{eq:truncatedhypotheses} Vh-V_\nu g\geq0,\qquad \kappa\in\Gamma_\ell\quad\text{when }\ell\geq2. \end{equation}
For $\ell=1$, assume static convexity alone. Then
\begin{equation}\label{eq:extendedadjacent} I_m^2\geq I_{m-1}I_{m+1},\qquad 0\leq m\leq\ell-1, \end{equation}
and, for integer orders $i<j<k$,
\begin{equation}\label{eq:interpolation} I_j^{\,k-i}\geq I_i^{\,k-j}I_k^{\,j-i},\qquad -1\leq i<j<k\leq\ell. \end{equation}
Equality in \eqref{eq:interpolation} holds if and only if $\Omega$ is a geodesic ball, whose center may differ from $p$.
\end{corollary}

\begin{proof}[Proof of \autoref{cor:interpolation}]
For $\ell=1$, both \eqref{eq:extendedadjacent} and \eqref{eq:interpolation} reduce to \eqref{eq:volume}, including its equality statement. Assume $\ell\geq2$.

For $m=0$, the assertion is \eqref{eq:volume}. For $m=1$, cone inclusion $\Gamma_\ell\subset\Gamma_2$ implies $T_1>0$, so $\mathcal Q\geq0$ by \eqref{eq:Qfactor}. \autoref{thm:allorders} therefore applies. For $2\leq m\leq\ell-1$, we have $\Gamma_\ell\subset\Gamma_{m+1}$ and apply \autoref{thm:allorders}. Since $I_{-1},I_0,\ldots,I_\ell$ are positive, the adjacent inequalities imply that their successive logarithmic differences are nonincreasing. Averaging over $i<m\leq j$ and $j<m\leq k$, we find
\[ \frac{\log I_j-\log I_i}{j-i}\geq\frac{\log I_k-\log I_j}{k-j}. \]
Rearranging and exponentiating proves \eqref{eq:interpolation}.

Equality in \eqref{eq:interpolation} forces equality at every adjacent order $i<m<k$. For each such $m\geq1$, the preceding argument verifies the hypotheses of \autoref{thm:pointwiseReilly}, including $T_{m-1}>0$ when $m\geq2$. \autoref{cor:hyperbolicrigidity} therefore applies if such an order occurs. The remaining case, $(i,j,k)=(-1,0,1)$, follows from the equality case of \eqref{eq:volume}. Conversely, the values in \eqref{eq:offcenterballcalibration} form a geometric progression in the curvature order, so every geodesic ball attains equality.
\end{proof}

\section{Static convexity and Brunn--Minkowski theory}\label{sec:staticdomains}

Throughout this section, the ambient space is $\Hhyp^{n+1}$ and $V=\cosh\operatorname{dist}(\cdot,p)$ for a fixed point $p$. We describe the projection geometry underlying the mixed-volume interpretation developed in \autoref{sec:projection}. Since spatial projection identifies weighted volume with Euclidean volume and static convexity with Euclidean convexity, the Brunn--Minkowski statements below follow from their classical Euclidean counterparts.

\begin{definition}\label{def:staticconvex}
A compact domain $\Omega$ with smooth connected boundary is static convex with respect to $V$ if $Vh-V_\nu g\geq0$ on its boundary. It is strictly static convex if this tensor is positive definite.
\end{definition}

The terminology is adapted from the convex static condition of Brendle--Wang \cite[Def.~3]{BrendleWang}, which additionally requires $V_\nu>0$. The tensor condition $Vh-V_\nu g\geq0$ was used in the weighted Minkowski inequality \cite{Xia}. The definition of static convexity in \cite{HuLi} also requires $V_\nu>0$.

In polar coordinates about $p$,
\begin{equation}\label{eq:spaceformstaticpotential} \bar g=\dd r^2+\sinh^2r\,g_{\mathbb S^n},\qquad V=\cosh r. \end{equation}
By \eqref{eq:hyperbolicstaticpotential}, $|V_\nu|<V$, so horospherical convexity $\kappa_i\geq1$ implies strict static convexity relative to every $p$.

The projection formulas below are due to Li--Xu \cite[Lems.~2.1--2.3]{LiXuIso}; see also \cite[Lem.~2.2 and Cor.~2.1]{HuLiWanXu}. Spatial projection in Minkowski spacetime was used earlier by Gibbons \cite[Sec.~5]{Gibbons}. We express these formulas in terms of a complete flat metric.

\begin{proposition}\label{prop:staticaffine}
On $\Hhyp^{n+1}$, define
\begin{equation}\label{eq:staticflatmetric} \widehat g=\bar g+\dd V\otimes\dd V. \end{equation}
This metric is complete and flat, with Levi-Civita connection and volume form given by
\begin{equation}\label{eq:staticconnection} \bar D_a b=\bar\nabla_a b+\bar g(a,b)\bar\nabla\log V, \end{equation}
\begin{equation}\label{eq:staticparallelvolume} \dd\mathrm{vol}_{\widehat g}=V\dd\mathrm{vol},\qquad \bar D(V\dd\mathrm{vol})=0. \end{equation}
Let $\widehat h$ be the second fundamental form of $\Sigma$ with respect to its outward $\widehat g$-unit normal. For vectors $a,b$ tangent to $\Sigma$,
\begin{equation}\label{eq:staticaffineshape} -\bar g(\bar D_a b,\nu)=h(a,b)-\frac{V_\nu}{V}g(a,b), \end{equation}
\begin{equation}\label{eq:staticflatshape} Vh-V_\nu g=\sqrt{V^2-V_\nu^2}\,\widehat h. \end{equation}
Thus static convexity is equivalent to $\widehat h\geq0$, and weighted volume agrees with the Riemannian volume of $\widehat g$.
\end{proposition}

\begin{proof}
The static identities imply $\operatorname{grad}_{\widehat g}V=V^{-2}\bar\nabla V$. By the Koszul formula,
\[ \widehat g(\bar D_a b-\bar\nabla_a b,c)=(\bar\nabla^2V)(a,b)\dd V(c)=V\bar g(a,b)\dd V(c). \]
This proves \eqref{eq:staticconnection}. With $\rho=\sinh r$,
\begin{equation}\label{eq:staticflatpolar} \widehat g=\cosh^2r\dd r^2+\sinh^2r\,g_{\mathbb S^n}=\dd\rho^2+\rho^2g_{\mathbb S^n}. \end{equation}
The chart \eqref{eq:canonicalprojection} extends smoothly across $p$ and is an isometry from $(\Hhyp^{n+1},\widehat g)$ onto $\mathbb R^{n+1}$. In the hyperboloid model, $\bar g=\sum_i\dd x_i^2-\dd V^2$ and $\widehat g=\sum_i\dd x_i^2$. Thus $\widehat g$ is complete and flat. The rank-one determinant formula gives $\dd\mathrm{vol}_{\widehat g}=\sqrt{1+|\bar\nabla V|^2}\dd\mathrm{vol}=V\dd\mathrm{vol}$, and its volume form is parallel.

Equation~\eqref{eq:staticaffineshape} follows from $\bar g(\bar\nabla_a b,\nu)=-h(a,b)$. The outward unit normal for $\widehat g$ is
\[ \widehat\nu=\frac{V}{\sqrt{V^2-V_\nu^2}}\left(\nu-\frac{V_\nu}{V^2}\bar\nabla V\right). \]
Since $\widehat g(\widehat\nu,\cdot)=V\bar g(\nu,\cdot)/\sqrt{V^2-V_\nu^2}$, contraction with $\bar D_a b$ proves \eqref{eq:staticflatshape}. The square-root factor is positive because $V^2-V_\nu^2\geq1$.
\end{proof}

For a compact domain with smooth connected boundary in $\mathbb R^{n+1}$, nonnegativity of the outward second fundamental form implies local convexity in graph coordinates. The induced boundary metric is complete by compactness, and comparison with an enclosing Euclidean sphere shows that the second fundamental form is positive definite at a contact point. Since $n\geq2$, van Heijenoort's theorem \cite{Heijenoort} applies and the domain is convex. The converse follows from supporting hyperplanes. Therefore
\begin{equation}\label{eq:projectionvolume} \int_\Omega V\dd\mathrm{vol}=|\pi(\Omega)|,\qquad Vh-V_\nu g\geq0\ \Longleftrightarrow\ \pi(\Omega)\text{ is convex}. \end{equation}
For the volume identity and the convexity equivalence, see \cite[Lem.~2.1]{LiXuIso} and \cite[Cor.~2.1]{HuLiWanXu}, respectively. The chart is determined by the metric and the chosen center up to an orthogonal transformation. The affine segments determined by $\bar D$ are
\begin{equation}\label{eq:staticaffinesegments} \gamma(t)=\pi^{-1}\bigl((1-t)\pi(X_0)+t\pi(X_1)\bigr),\qquad 0\leq t\leq1. \end{equation}

For two static convex domains, define their static affine interpolation by
\begin{equation}\label{eq:staticaddition} \pi(\Omega_\lambda^{\mathrm{stat}})=(1-\lambda)\pi(\Omega_0)+\lambda\pi(\Omega_1),\qquad 0\leq\lambda\leq1. \end{equation}
The classical Brunn--Minkowski theorem \cite{Gardner} takes the following weighted form under spatial projection.

\begin{theorem}\label{thm:weightedBM}
Let $\Omega_0,\Omega_1\subset\Hhyp^{n+1}$ be compact domains with smooth connected strictly static convex boundaries, with respect to the same weight $V=\cosh\operatorname{dist}(\cdot,p)$. Then every $\Omega_\lambda^{\mathrm{stat}}$ has smooth connected strictly static convex boundary, and
\begin{equation}\label{eq:staticBM} I_{-1}(\Omega_\lambda^{\mathrm{stat}})^{1/(n+1)}\geq(1-\lambda)I_{-1}(\Omega_0)^{1/(n+1)}+\lambda I_{-1}(\Omega_1)^{1/(n+1)}. \end{equation}
For $0<\lambda<1$, equality holds if and only if $\pi(\Omega_1)=a\pi(\Omega_0)+z$ for some $a>0$ and $z\in\mathbb R^{n+1}$.
\end{theorem}

For the mixed-volume formulation, recall that a convex body is a compact convex set with nonempty interior. We denote Euclidean mixed volume by $\operatorname{MV}_E$, normalized by $\operatorname{MV}_E(K[n+1])=|K|$.

For the interpolation \eqref{eq:staticaddition}, define weighted mixed volume by
\begin{equation}\label{eq:staticmixedvolume} \mathcal V(\Omega_0,\ldots,\Omega_n)=(n+1)\operatorname{MV}_E\bigl(\pi(\Omega_0),\ldots,\pi(\Omega_n)\bigr). \end{equation}
This functional is symmetric and multilinear with respect to Minkowski addition and nonnegative scaling of the projected bodies, with $\mathcal V(\Omega[n+1])=I_{-1}(\Omega)$. The classical Brunn--Minkowski and Alexandrov--Fenchel inequalities \cite{Gardner,ShenfeldVanHandel} take the following form in these coordinates.
\begin{proposition}\label{prop:staticBM}
Let the domains below be compact, with smooth connected strictly static convex boundaries and the same weight center $p$. The following statements hold.
\begin{enumerate}[(i)]
\item Let $0\leq k\leq n-1$ and fix $\Lambda_1,\ldots,\Lambda_k$. The functional
\begin{equation}\label{eq:staticmixedBM} \Omega\longmapsto\mathcal V(\Omega[n+1-k],\Lambda_1,\ldots,\Lambda_k)^{1/(n+1-k)} \end{equation}
is concave along static affine segments. When $k=0$, the list of fixed arguments is empty and this is \eqref{eq:staticBM}.
\item With a fixed list of $n-1$ domains in the remaining arguments, the mixed-volume inequality is
\begin{equation}\label{eq:staticAF} \mathcal V(\Omega_0,\Omega_1,\ldots)^2\geq\mathcal V(\Omega_0,\Omega_0,\ldots)\mathcal V(\Omega_1,\Omega_1,\ldots). \end{equation}
\end{enumerate}
\end{proposition}

\begin{proof}[Proof of \autoref{thm:weightedBM} and \autoref{prop:staticBM}]
The volume identity \eqref{eq:projectionvolume} reduces \eqref{eq:staticBM} to the classical Brunn--Minkowski theorem \cite{Gardner}, which also gives its homothety equality case. For a smooth function $v$ on $\mathbb S^n$, write $\tau[v]=\nabla_{\mathbb S}^2v+v\mathrm{Id}$. For smooth strictly static convex endpoints, their support functions satisfy $\tau[u_i]>0$. The interpolating support function has $\tau[u_\lambda]=(1-\lambda)\tau[u_0]+\lambda \tau[u_1]>0$, so its lift is smooth and strictly static convex by \autoref{prop:staticaffine}.

Under \eqref{eq:staticmixedvolume}, \eqref{eq:staticAF} is the Euclidean Alexandrov--Fenchel inequality \cite[Thm.~1.1]{ShenfeldVanHandel}, and \eqref{eq:staticmixedBM} is the generalized Brunn--Minkowski inequality for mixed volumes \cite{Schneider}. The latter applies because the projected bodies have nonempty interiors and static affine interpolation becomes ordinary Minkowski interpolation under $\pi$.
\end{proof}

The affine supporting hyperplanes are $x\cdot\theta=c$, where $\theta\in\mathbb S^n$. Set $\varphi_\theta(X)=x\cdot\theta$. In the hyperboloid model,
\[ \bar\nabla^2\varphi_\theta=\varphi_\theta\bar g,\qquad |\bar\nabla\varphi_\theta|^2=1+\varphi_\theta^2. \]
Moreover, $\bar g(\bar\nabla V,\bar\nabla\varphi_\theta)=V\varphi_\theta$. Thus $\operatorname{arsinh}\varphi_\theta$ is signed distance to the totally geodesic hyperplane $\varphi_\theta=0$, and on the level set $\varphi_\theta=c$, with normal $\bar\nabla\varphi_\theta/\sqrt{1+c^2}$,
\begin{equation}\label{eq:staticaffinehyperplanes} h=\frac{c}{\sqrt{1+c^2}}g,\qquad \frac{V_\nu}{V}=\frac{c}{\sqrt{1+c^2}},\qquad Vh-V_\nu g=0. \end{equation}
These equidistant hypersurfaces serve as affine supporting hypersurfaces for static convex domains.

Let $K=\pi(\Omega)$ be the convex spatial projection of a static convex domain. Its support function is
\[ u(\theta)=\sup_{x\in K}x\cdot\theta,\qquad \theta\in\mathbb S^n. \]
At a boundary point whose outward Euclidean normal is $\theta$, the normal formula \cite[Lem.~2.2]{LiXuIso} implies
\begin{equation}\label{eq:staticdomainnormal} \frac{V_\nu}{V}=\frac{u(\theta)}{\sqrt{1+u(\theta)^2}}. \end{equation}
Formula~\eqref{eq:staticdomainnormal} remains valid for $u<0$. If $p\in\operatorname{int}(\Omega)$, then $u>0$ and $V_\nu>0$, so static convexity implies $h\geq(V_\nu/V)g>0$. Since $K$ is convex and $\pi$ preserves radial rays, $\Omega$ is also star-shaped with respect to $p$. Thus static convexity and $p\in\operatorname{int}(\Omega)$ imply the positivity assumptions in \cite{HuLi}.

\section{Mixed-volume interpretation and the Hu--Li comparisons}\label{sec:projection}

Throughout this section, $\Omega\subset\Hhyp^{n+1}$ is a compact domain with smooth connected boundary $\Sigma$, and $V=\cosh\operatorname{dist}(\cdot,p)$. Under $h>0$ and $S>0$, we interpret the adjacent integral inequalities \eqref{eq:allorders} and the weighted-volume inequality \eqref{eq:volume} proved above as Euclidean mixed-volume inequalities. Both the integral inequalities \eqref{eq:BMcurvature} and the representation \eqref{eq:originalmixedvolumes} extend to $h\geq0$ and $S\geq0$, by outer parallel approximation and \autoref{lem:weakmixed}, respectively. As a geometric application, the representation supplies the highest-order estimate \eqref{eq:HuLiterminal}, which, combined with log-concavity, recovers the comparisons in \autoref{thm:HuLiconjecture}. For related weighted inequalities involving convex weights, volume, and quermassintegrals, see Kwong--Wei \cite{KwongWei} and Wu \cite{Wu}.

\subsection{The normal body and weak convexity}\label{subsec:mixedcurvature}

The normal-body representation uses the projection formulas in \cite[Lems.~2.1--2.3]{LiXuIso}.

\begin{theorem}\label{thm:mixedidentification}
Let $\Omega\subset\Hhyp^{n+1}$ be a compact domain with smooth connected boundary, and let $V=\cosh\operatorname{dist}(\cdot,p)$. Assume $h>0$ and $Vh-V_\nu g>0$. Let $K=\pi(\Omega)$ and let $u=u_K$ be its Euclidean support function. Then
\begin{equation}\label{eq:normalbodysupport} u_L=\sqrt{1+u^2} \end{equation}
is the support function of a smooth convex body $L$ containing the Euclidean unit ball, and
\begin{equation}\label{eq:originalmixedvolumes} I_j=(n+1)\operatorname{MV}_E(K[n-j],L[j+1]),\qquad -1\leq j\leq n. \end{equation}
Both bodies have smooth boundaries with positive principal curvatures, and all $I_j$ are positive. The point $p$ may lie outside $\Omega$.
\end{theorem}

The assumption $h>0$ ensures convexity of the normal body, as the proof below shows, and \autoref{lem:weakmixed} extends the representation to $h\geq0$, $S\geq0$.

\begin{proof}[Proof of \autoref{thm:mixedidentification}]
Under the assumptions $h>0$ and $Vh-V_\nu g>0$, \autoref{prop:staticaffine} implies that $K$ is a smooth convex body with positive Euclidean principal curvatures and support parametrization $x(\theta)=\nabla_{\mathbb S}u+u\theta$. Identifying $T_\theta\mathbb S^n$ with $\theta^\perp$, set
\[ \tau[u]=\mathrm{d}x=\nabla_{\mathbb S}^2u+u\mathrm{Id}>0. \]
At the lifted point $X=(V,x)$, the normal formula \cite[Lem.~2.2]{LiXuIso} is
\begin{equation}\label{eq:spaceformnormal} \nu=\frac{(uV,ux+\theta)}{u_L},\qquad \frac{V_\nu}{V}=\frac{u}{u_L},\qquad u_L=\sqrt{1+u^2}. \end{equation}
Thus the spatial normal image satisfies
\[ y=\pi\nu=\nabla_{\mathbb S}u_L+u_L\theta,\qquad \tau[u_L]=\mathrm{d}y=\nabla_{\mathbb S}^2u_L+u_L\mathrm{Id}. \]
Here $\pi\nu$ denotes the spatial component of the ambient normal. The restriction $\mathrm{d}\pi|_{T_X\Sigma}$ maps $T_X\Sigma$ isomorphically onto $\theta^\perp$. By the ambient Weingarten formula $\mathrm{d}\nu=h\,\mathrm{d}X$,
\begin{equation}\label{eq:normalprojection} \tau[u_L]\tau[u]^{-1}=\bigl(\mathrm{d}\pi|_{T_X\Sigma}\bigr)h\bigl(\mathrm{d}\pi|_{T_X\Sigma}\bigr)^{-1}. \end{equation}
The algebraic identity \eqref{eq:normalprojection} remains valid under strict static convexity alone; $h>0$ is used below to deduce $\tau[u_L]>0$.
The left-hand side is similar to the symmetric matrix $\tau[u]^{-1/2}\tau[u_L]\tau[u]^{-1/2}$, whose eigenvalues are positive because $h>0$, so $\tau[u_L]>0$. By the support function criterion \cite[Lem.~2.1]{ShenfeldVanHandel}, $u_L$ is the support function of a smooth convex body $L$ with positive principal curvatures. Since $u_L\geq1$, this body contains the unit ball.

Let $\dd\mu_{\mathbb S}$ be round measure on $\mathbb S^n$. The boundary determinant formula for $\widehat g$ and \eqref{eq:spaceformnormal} imply
\[ \dd\widehat\mu=\det \tau[u]\dd\mu_{\mathbb S}=\sqrt{V^2-V_\nu^2}\dd\mu=\frac{V}{u_L}\dd\mu. \]
This is the area formula in \cite[Lem.~2.3]{LiXuIso}, also recorded in \cite[Lem.~2.2]{HuLiWanXu}. Consequently,
\begin{equation}\label{eq:weightedGaussmeasure} V\dd\mu=u_L\det \tau[u]\dd\mu_{\mathbb S}. \end{equation}
Taking determinants in \eqref{eq:normalprojection} and comparing coefficients, we obtain
\begin{equation}\label{eq:mixedGausscurvature} D\bigl(\tau[u][n-j],\tau[u_L][j]\bigr)=E_j\det \tau[u],\qquad 0\leq j\leq n. \end{equation}
Hence
\[ I_j=\int_{\mathbb S^n}u_L D\bigl(\tau[u][n-j],\tau[u_L][j]\bigr)\dd\mu_{\mathbb S}. \]
The mixed-volume integral formula \cite[Eq.~(1.4)]{ShenfeldVanHandel} gives \eqref{eq:originalmixedvolumes} for $0\leq j\leq n$. The case $j=-1$ follows from \eqref{eq:projectionvolume}. Positivity follows because both bodies have nonempty interior.
\end{proof}

The following approximation extends the representation to the weak hypotheses used in the Hu--Li comparisons.

\begin{lemma}\label{lem:weakmixed}
Let $\Omega\subset\Hhyp^{n+1}$ be a compact domain with smooth connected boundary, let $V=\cosh\operatorname{dist}(\cdot,p)$, and assume $h\geq0$ and $Vh-V_\nu g\geq0$. Set $K=\pi(\Omega)$ and $u=u_K$. There are smooth domains $\Omega_m$ with $h_m>0$ and $Vh_m-V_{\nu_m}g_m>0$, for the same potential $V$, whose boundaries converge to $\partial\Omega$ in $C^2$. Their projected bodies and normal bodies converge in Hausdorff distance to $K$ and a convex body $L$ with support function $u_L=\sqrt{1+u^2}$. The mixed-volume representation \eqref{eq:originalmixedvolumes} holds in the limit, and $I_j>0$ for $-1\leq j\leq n$.
\end{lemma}

\begin{proof}
The Klein chart maps geodesics to straight segments and preserves the sign of the second fundamental form. Thus $h\geq0$ makes the Klein image of $\Sigma$ locally convex. This image is compact and connected, its induced Euclidean metric is complete, and comparison with an enclosing Euclidean sphere supplies a strictly convex point. Van Heijenoort's theorem \cite{Heijenoort} therefore implies that it bounds a convex body, so $\Omega$ is geodesically convex. For $t>0$, the map $X\mapsto(\cosh t)X+(\sinh t)\nu$ parametrizes its outer parallel boundary. This map is injective by uniqueness of the nearest point to a closed convex set, and its tangential differential is nonsingular because $\cosh t+\kappa_i\sinh t>0$ for every $i$. Thus the parallel domains $\Omega_t$ are smooth and converge smoothly to $\Omega$ as $t\downarrow0$. Along their normal geodesics,
\begin{equation}\label{eq:HuLiweakparallel} \kappa_i(t)=\frac{\kappa_i+\tanh t}{1+\kappa_i\tanh t}>0,\qquad V(t)\kappa_i(t)-V_\nu(t)=\frac{V\kappa_i-V_\nu}{\cosh t+\kappa_i\sinh t}\geq0. \end{equation}
The outer parallels satisfy $h_t>0$, but the second formula preserves zeros of the static convexity tensor, so we next approximate their spatial projections to obtain strict static convexity.

Fix $t>0$ and translate $K_t=\pi(\Omega_t)$ so that the origin lies in its interior. Let $F$ be the Minkowski functional of the translated body. Radial rays meet its smooth boundary transversely, so $F$ is smooth away from the origin. Its square is convex, and $G_\eta=F^2+\eta|x|^2$ satisfies $\nabla_E^2G_\eta\geq2\eta\mathrm{Id}$ away from the origin for $\eta>0$, where $\nabla_E$ is the Euclidean connection. The unit level has the radial parametrization
\[ x_\eta(\theta)=\frac{\theta}{\sqrt{F(\theta)^2+\eta}},\qquad \theta\in\mathbb S^n. \]
Since $x\cdot\nabla_EG_\eta=2$ on this level, the unit sublevel set has smooth boundary with positive Euclidean principal curvatures. The parametrizations $x_\eta$ converge to $x_0$ in $C^2$ as $\eta\downarrow0$. Translate back and lift by $\pi^{-1}$. The lifts are strictly static convex by \autoref{prop:staticaffine}. Since $\kappa_i(t)\geq\tanh t>0$ by \eqref{eq:HuLiweakparallel}, they retain $h>0$ for sufficiently small $\eta$ by $C^2$ continuity. Choose $t_m\downarrow0$ and then $\eta_m$ so that $h>0$ and the lifted boundary is within $1/m$ of $\partial\Omega_{t_m}$ in $C^2$. This produces the stated sequence.

Let $K_m=\pi(\Omega_m)$ and let $u_m=u_{K_m}$ be their support functions. Convergence in $C^2$ of their boundaries implies Hausdorff convergence and uniform convergence $u_m\to u$. The preceding construction yields normal bodies $L_m$ with support functions $u_{L_m}=\sqrt{1+u_m^2}$. Since $u_{L_m}\geq1$, the corresponding convex bodies contain the unit ball. Uniform convergence $u_{L_m}\to\sqrt{1+u^2}$ therefore defines a convex body $L$ containing the unit ball. Mixed volumes are continuous under Hausdorff convergence \cite[Sec.~5.1]{Schneider}, whereas the boundary integrals $I_j$, $0\leq j\leq n$, are continuous under $C^2$ boundary convergence. The weighted volumes converge by \eqref{eq:projectionvolume} and $|K_m|\to|K|$. Passing to the limit proves \eqref{eq:originalmixedvolumes}.

The body $K$ has nonempty interior because $\pi$ is a diffeomorphism, and $L$ contains the unit ball. Translation invariance and monotonicity of mixed volumes \cite[Sec.~5.1]{Schneider} therefore imply $I_j>0$ for every $-1\leq j\leq n$.
\end{proof}

By the Euclidean volume polynomial, we have
\begin{equation}\label{eq:originalvolumepolynomial} (n+1)|\alpha K+\beta L|=\sum_{\ell=0}^{n+1}\binom{n+1}{\ell}I_{\ell-1}\alpha^{n+1-\ell}\beta^\ell,\qquad \alpha,\beta\geq0. \end{equation}

\subsection{Log-concavity and its mixed-volume interpretation}\label{subsec:mixedconcavity}

\begin{corollary}\label{cor:mixedpolynomial}
Under the hypotheses of \autoref{thm:mixedidentification}, the following statements hold.
\begin{enumerate}[(i)]
\item The sequence $I_{-1},I_0,\ldots,I_n$ is positive and log-concave:
\begin{equation}\label{eq:BMcurvature} I_j^2\geq I_{j-1}I_{j+1},\qquad 0\leq j\leq n-1. \end{equation}
Consequently, for integer orders $-1\leq k<\ell<m\leq n$,
\begin{equation}\label{eq:BMinterpolation} I_\ell^{\,m-k}\geq I_k^{\,m-\ell}I_m^{\,\ell-k}. \end{equation}
\item Equality at any order in \eqref{eq:BMcurvature}, or for any triple in \eqref{eq:BMinterpolation}, holds if and only if $\Omega$ is a geodesic ball, whose center may differ from $p$.
\item The function
\begin{equation}\label{eq:polynomialBM} (\alpha,\beta)\longmapsto\left[\sum_{\ell=0}^{n+1}\binom{n+1}{\ell}I_{\ell-1}\alpha^{n+1-\ell}\beta^\ell\right]^{1/(n+1)} \end{equation}
is concave on the nonnegative quadrant.
\end{enumerate}
Inequalities \eqref{eq:BMcurvature} and \eqref{eq:BMinterpolation}, together with the concavity in \eqref{eq:polynomialBM}, extend by approximation to $h\geq0$ and $S\geq0$. The equality statement in \textup{(ii)} assumes $h>0$ and $S>0$.
\end{corollary}

\begin{proof}
Assume first $h\geq0$ and $S\geq0$. The outer parallel domains $\Omega_t$ have smooth connected boundaries and satisfy $h_t>0$ and $S_t\geq0$ for $t>0$ by \eqref{eq:HuLiweakparallel}. Applying \autoref{thm:allorders} for $1\leq j\leq n-1$ and \eqref{eq:volume} for $j=0$ to $\Omega_t$, then letting $t\downarrow0$, proves \eqref{eq:BMcurvature} by smooth convergence. The volume integral $I_{-1}$ is positive, as are the boundary integrals: $h\geq0$, and \eqref{eq:canonicalellipticpoint} gives $h>0$ near a farthest point from $p$. Thus the successive logarithmic differences are nonincreasing, which implies \eqref{eq:BMinterpolation}. Equality in \eqref{eq:BMinterpolation} forces equality at an intermediate adjacent order.

For the convex bodies $K,L$ supplied by \autoref{lem:weakmixed}, the representation \eqref{eq:originalmixedvolumes} also gives a geometric proof of \eqref{eq:BMcurvature}: apply Alexandrov--Fenchel with remaining arguments $K[n-j-1],L[j]$. Concavity of \eqref{eq:polynomialBM} follows from \eqref{eq:originalvolumepolynomial} and the Euclidean Brunn--Minkowski inequality.

Under the strict hypotheses, $K,L$ and every body in the fixed list $K[n-j-1],L[j]$ have smooth boundaries with positive principal curvatures. The smooth equality theorem for Alexandrov--Fenchel \cite[Sec.~6.2]{ShenfeldVanHandel} therefore applies at every $0\leq j\leq n-1$: equality holds if and only if $L=aK+z$ for some $a>0$ and $z\in\mathbb R^{n+1}$. Their support functions satisfy $u_L=au+z\cdot\theta$, so $\tau[u_L]=a\tau[u]$ and $h=a\mathrm{Id}$ by \eqref{eq:normalprojection}. A compact connected totally umbilical hypersurface in hyperbolic space is a geodesic sphere. Conversely, \eqref{eq:offcenterballcalibration} shows that every geodesic ball attains equality at every adjacent order and for every triple.
\end{proof}

The cases $j=0$ and $j=1$ in \eqref{eq:BMcurvature} are \eqref{eq:volume} and \eqref{eq:three}, respectively. To derive the adjacent inequalities directly from Brunn--Minkowski concavity, assume $h>0$ and strict static convexity and set $\Omega^\nu=\pi^{-1}(L)$. Then
\begin{equation}\label{eq:staticnormalmixed} I_j(\Omega)=\mathcal V\bigl(\Omega[n-j],\Omega^\nu[j+1]\bigr),\qquad -1\leq j\leq n. \end{equation}
The normal domain is smooth and strictly static convex, contains $p$, and remains fixed in the following variation.

\begin{proposition}\label{prop:normalrayBM}
Let $\Omega$ be a compact domain with smooth connected boundary, and assume $h>0$ and $Vh-V_\nu g>0$. The domains $\Omega_s$ defined by $\pi(\Omega_s)=K+sL$, $s\geq0$, have smooth strictly static convex boundaries. For $0\leq k\leq n-1$,
\begin{equation}\label{eq:adjacentquadraticBM} \mathcal V\bigl(\Omega_s[2],\Omega[n-k-1],\Omega^\nu[k]\bigr)=I_{k-1}+2sI_k+s^2I_{k+1}. \end{equation}
The square root of this polynomial is concave, and
\begin{equation}\label{eq:adjacentsqrtdeficit} \frac{d^2}{ds^2}\sqrt{I_{k-1}+2sI_k+s^2I_{k+1}}=\frac{I_{k-1}I_{k+1}-I_k^2}{(I_{k-1}+2sI_k+s^2I_{k+1})^{3/2}}\leq0. \end{equation}
For any fixed order $k$, the second derivative in \eqref{eq:adjacentsqrtdeficit} vanishes at some $s\geq0$ if and only if $\Omega$ is a geodesic ball, whose center may differ from $p$.
\end{proposition}

\begin{proof}
The support function of $K+sL$ is $u+s u_L$, and $\tau[u+s u_L]=\tau[u]+s\tau[u_L]>0$. Thus every $\Omega_s$ has smooth strictly static convex boundary by \autoref{prop:staticaffine}.

Multilinearity and \eqref{eq:staticnormalmixed} give \eqref{eq:adjacentquadraticBM}. On each finite interval of $s$, the domains form a static affine segment, so \autoref{prop:staticBM} with the fixed arguments $\Omega[n-k-1],\Omega^\nu[k]$ proves concavity of the square root. Differentiation yields \eqref{eq:adjacentsqrtdeficit}, whose equality case follows from \autoref{cor:mixedpolynomial}.
\end{proof}

The second derivative in \eqref{eq:adjacentsqrtdeficit} is a negative multiple of $I_k^2-I_{k-1}I_{k+1}$. For $k\geq1$, \eqref{eq:exactdeficit} decomposes this deficit divided by $I_{k-1}$; at $k=0$, \eqref{eq:XiaReillycomparison} decomposes a positive multiple of $I_0^2-I_{-1}I_1$. The normal body $L$ and the remaining domains in \eqref{eq:adjacentquadraticBM} are fixed, so this mixed polynomial generally differs from $I_{k-1}(\Omega_s)$.

\subsection{The Hu--Li comparisons}\label{subsec:HuLicomparisons}

The proof of \autoref{thm:HuLiconjecture} combines the log-concavity inequalities \eqref{eq:BMcurvature} with the highest-order estimate \eqref{eq:HuLiterminal}. The $j$th weighted curvature integral in \cite{HuLi} is $I_{j-1}$ in our notation.

\begin{proof}[Proof of \autoref{thm:HuLiconjecture}]
Let $K,L$ be the convex bodies supplied by \autoref{lem:weakmixed}, and let $S_L$ be the Euclidean surface area measure of $L$, with $|\partial L|=S_L(\mathbb S^n)$. The mixed-volume integral formula for convex bodies \cite[Sec.~5.1]{Schneider}, together with \eqref{eq:originalmixedvolumes}, gives
\[ I_{n-1}=\int_{\mathbb S^n}u\dd S_L,\qquad I_n=\int_{\mathbb S^n}\sqrt{1+u^2}\dd S_L=(n+1)|L|. \]
By Jensen's inequality for $s\mapsto\sqrt{1+s^2}$ on $\mathbb R$ with respect to the probability measure $S_L/|\partial L|$, as in \cite[Eqs.~(2.17)--(2.18)]{LiXuIso}, we have $I_n^2-I_{n-1}^2\geq|\partial L|^2$. The Euclidean isoperimetric inequality for convex bodies \cite[Sec.~7.2]{Schneider} then implies
\begin{equation}\label{eq:HuLiterminal} I_n^2-I_{n-1}^2\geq|\partial L|^2\geq\omega_n^{\frac{2}{n+1}}I_n^{\frac{2n}{n+1}}. \end{equation}
If the first and last terms in \eqref{eq:HuLiterminal} are equal, both inequalities are equalities. Since $L$ has nonempty interior, the isoperimetric equality case makes it a Euclidean ball. Strict convexity of $s\mapsto\sqrt{1+s^2}$ makes $u$ constant $S_L$-almost everywhere. The measure $S_L$ is then a positive multiple of round measure, so continuity implies $u\equiv c>0$, with positivity following from $I_{n-1}>0$. Thus $K=B_c(0)$ and $\Omega=B_{\operatorname{arsinh}c}(p)$.

For $-1\leq j\leq n$, define $R_j>0$ by $I_j=I_j(B_{R_j}(p))$. The centered-ball values are
\begin{equation}\label{eq:centeredballcalibration} I_j(B_R(p))=\omega_n\sinh^{n-j} R\cosh^{j+1} R,\qquad -1\leq j\leq n. \end{equation}
Their logarithmic derivatives are $(n-j)\coth R+(j+1)\tanh R>0$. For $j<n$ the range is $(0,\infty)$, while the range at $j=n$ is $(\omega_n,\infty)$. The integrals are positive for $j<n$ by \autoref{lem:weakmixed}. At the top order, \eqref{eq:HuLiterminal} and $I_{n-1}>0$ imply $I_n>\omega_n$, so $R_n$ also exists uniquely. Substituting $I_n=\omega_n\cosh^{n+1} R_n$ into \eqref{eq:HuLiterminal} shows $I_{n-1}\leq I_{n-1}(B_{R_n}(p))$, and therefore $R_{n-1}\leq R_n$.

Suppose $R_j\leq R_{j+1}$ for some $0\leq j\leq n-1$. The adjacent integral inequality \eqref{eq:BMcurvature}, valid under the weak hypotheses by \autoref{cor:mixedpolynomial}, together with monotonicity of the ball values and their geometric progression in the curvature order, implies
\begin{equation}\label{eq:HuLirecursion} I_{j-1}\leq\frac{I_j^2}{I_{j+1}}\leq\frac{I_j(B_{R_j}(p))^2}{I_{j+1}(B_{R_j}(p))}=I_{j-1}(B_{R_j}(p)). \end{equation}
Thus $R_{j-1}\leq R_j$. Backward induction, including $j=0$, proves $R_{-1}\leq R_0\leq\cdots\leq R_n$, which is \eqref{eq:HuLiconjecture}.

If equality holds for $k<\ell$, then $R_k=R_\ell$ and hence $R_k=R_{k+1}$. Whenever $R_{j-1}=R_j$, with $0\leq j\leq n-1$, we have
\[
I_{j+1}\leq\frac{I_j^2}{I_{j-1}}=I_{j+1}(B_{R_j}(p))\leq I_{j+1}.
\]
Here the first inequality is \eqref{eq:BMcurvature}, the equality follows from \eqref{eq:centeredballcalibration} and $R_{j-1}=R_j$, and the last inequality uses $R_j\leq R_{j+1}$. Thus $R_j=R_{j+1}$. Starting from $R_k=R_{k+1}$, we obtain $R_{n-1}=R_n$, including when $k=-1$. By \eqref{eq:centeredballcalibration}, the first and last terms in \eqref{eq:HuLiterminal} are then equal. Its equality case implies that $\Omega$ is a ball centered at $p$. Conversely, every such ball attains equality.
\end{proof}

\begin{remark}\label{rem:HuLihypotheses}
The assumptions in \eqref{eq:HuLihypotheses} are sufficient for all the inequalities in \autoref{thm:HuLiconjecture}, but need not be necessary for a particular pair. For $(k,\ell)=(-1,0)$, Li--Xu \cite[Cor.~1.1]{LiXuIso} require only a bounded domain with $C^1$ boundary. For $k=-1$ and $0\leq j\leq n$, \autoref{thm:HuLiconjecture} reads
\begin{equation}\label{eq:HuLivolume} I_j\geq\omega_n\left(\frac{I_{-1}}{\omega_n}\right)^{\frac{n-j}{n+1}}\left[1+\left(\frac{I_{-1}}{\omega_n}\right)^{\frac{2}{n+1}}\right]^{\frac{j+1}{2}}. \end{equation}
On the original static convex class, this is the weighted-volume comparison in \cite[Thm.~1.7]{HuLi}. The present statement also includes semidefinite tensors in \eqref{eq:HuLihypotheses} and weight centers outside the domain.
\end{remark}

\section{Support variations and parallel domains}\label{sec:variations}\label{subsec:variations}

Throughout this section, $\Omega\subset\Hhyp^{n+1}$ is a compact domain with smooth connected boundary $\Sigma$, and $V=\cosh\operatorname{dist}(\cdot,p)$. We consider variations of the spatial support function and outer parallel domains, relating the deficit $I_0^2-I_{-1}I_1$ in the weighted-volume inequality to the second variations along these two families.

\subsection{Static support variations}

If $\partial_tX=f\nu$, then $I_{-1}'=(n+1)\int_{\Sigma_t}Vf\dd\mu_t$. The curvature variation formula in \cite[Prop.~3.1]{HuLi} reads
\begin{equation}\label{eq:intrinsicvariation} \frac{d}{dt}I_j(\Omega_t)=\int_{\Sigma_t}\bigl[(n-j)VE_{j+1}+(j+1)V_\nu E_j\bigr]f\dd\mu_t. \end{equation}
This holds for $0\leq j\leq n$, with the first term omitted at $j=n$.

The quadratic inequality \eqref{eq:staticquadraticinequality} below expresses the local Minkowski inequality in terms of hyperbolic normal speeds.

\begin{proposition}\label{prop:staticsecondvariation}
Let $\Omega$ be a compact domain with smooth connected boundary satisfying $Vh-V_\nu g>0$, and let $u$ be the Euclidean support function of $K=\pi(\Omega)$. For every smooth function $f$ on $\Sigma$,
\begin{equation}\label{eq:staticquadraticinequality} \int_\Sigma\bigl[nVE_1f^2-V^2(S^{-1})^{ij}\nabla_i f\nabla_j f\bigr]\dd\mu\leq\frac{n}{I_{-1}}\left(\int_\Sigma Vf\dd\mu\right)^2. \end{equation}
In the Gauss parametrization by the outward Euclidean normal $\theta\in\mathbb S^n$ of $\partial K$, equality holds precisely for
\begin{equation}\label{eq:staticquadraticnull} f(\theta)=\frac{a u(\theta)+z\cdot\theta}{\sqrt{1+u(\theta)^2}},\qquad a\in\mathbb R,\quad z\in\mathbb R^{n+1}. \end{equation}
\end{proposition}

\begin{proof}
Put $\widetilde u=\sqrt{1+u^2}$ and regard $v=\widetilde u f$ as a function on $\mathbb S^n$. The functions $u_s=u+sv$ have a positive definite curvature-radius matrix for sufficiently small $|s|$, and their lifts retain strict static convexity. At fixed Euclidean Gauss direction, their normal speeds satisfy
\[ f_s=\frac{v}{\sqrt{1+u_s^2}},\qquad \partial_sf_s=-\frac{V_\nu}{V}f_s^2. \]
All geometric quantities are evaluated on the moving boundary. To compute the tangential velocity, consider the function $\varphi_\theta(X)=x\cdot\theta$ on the hyperboloid. It satisfies $\bar\nabla^2\varphi_\theta=\varphi_\theta\bar g$ and $|\bar\nabla\varphi_\theta|^2=1+\varphi_\theta^2$. Hence the signed distance $b_\theta=\operatorname{arsinh}\varphi_\theta$ to the corresponding totally geodesic hyperplane satisfies
\[ \bar\nabla^2b_\theta=\tanh b_\theta\bigl(\bar g-\dd b_\theta\otimes\dd b_\theta\bigr). \]
At the supporting point, $\bar\nabla b_\theta=\nu$ and $\tanh b_\theta=u_s/\sqrt{1+u_s^2}=V_\nu/V$. Write $\partial_sX=f_s\nu+W_s$. Differentiating the normal along this variation, we have $(V_\nu/V)W_s=-\nabla f_s+hW_s$. Therefore,
\[ W_s=VS^{-1}\nabla f_s. \]
The derivative $\partial_sf_s$ above is taken at fixed Euclidean Gauss direction, so the tangential velocity must also be included when differentiating the integral. By the transport formula and integration of $\diver(Vf_sW_s)$, we have
\[
\frac1{n+1}\frac{d^2}{ds^2}I_{-1}(\Omega_s)=\int_{\Sigma_s}\bigl[V\partial_sf_s+(V_\nu+nVE_1)f_s^2-V\langle W_s,\nabla f_s\rangle\bigr]\dd\mu_s.
\]
Substituting $\partial_sf_s=-(V_\nu/V)f_s^2$ cancels the $V_\nu$ term. With $W_s=VS^{-1}\nabla f_s$, we obtain
\begin{equation}\label{eq:staticsecondvariation} \frac{1}{n+1}\frac{d^2}{ds^2}I_{-1}(\Omega_s)=\int_{\Sigma_s}\bigl[nVE_1f_s^2-V^2(S^{-1})^{ij}\nabla_i f_s\nabla_j f_s\bigr]\dd\mu_s. \end{equation}
By the weighted-volume Brunn--Minkowski inequality \eqref{eq:staticBM}, $I_{-1}(\Omega_s)^{1/(n+1)}$ is concave. Differentiating twice at $s=0$ and using $I_{-1}'=(n+1)\int Vf_s\dd\mu$, we conclude \eqref{eq:staticquadraticinequality}.

To justify the equality case for an arbitrary smooth support variation, let $A_v(s)=(n+1)|K_s|$, where $K_s$ has support $u+sv$, and write $a_j=A_v^{(j)}(0)$ for $j=0,1,2$. Homogeneity of volume shows that replacing $v$ by $v+cu$ changes these derivatives to $a_0$, $a_1+c(n+1)a_0$, and $a_2+2cna_1+c^2n(n+1)a_0$. Consequently, the Brunn--Minkowski second-variation deficit $na_1^2-(n+1)a_0a_2$ is unchanged. Choose $c$ large enough that $\tau[v]+c\tau[u]>0$, so that $v+cu$ is the support function of a smooth convex body with positive principal curvatures. Smooth mixed-volume rigidity \cite[Sec.~6.2]{ShenfeldVanHandel} then shows that equality holds if and only if $v+cu$ differs from a positive multiple of $u$ by a linear function. Equivalently, $v=au+z\cdot\theta$ with $a\in\mathbb R$ and $z\in\mathbb R^{n+1}$, as asserted in \eqref{eq:staticquadraticnull}.
\end{proof}

Setting $f=1$ in \eqref{eq:staticquadraticinequality}, we recover the weighted-volume inequality \eqref{eq:volume}. Equality is equivalent to $\widetilde u=au+z\cdot\theta$. Integrating this identity over $\mathbb S^n$ shows that $a>0$, since $\widetilde u>0$ and $\int_{\mathbb S^n}u\dd\mu_{\mathbb S}=n^{-1}\int_{\mathbb S^n}\tr \tau[u]\dd\mu_{\mathbb S}>0$. Thus $\tau[\widetilde u]=a\tau[u]$ and the algebraic identity \eqref{eq:normalprojection} implies $h=a\mathrm{Id}$. A compact connected totally umbilical hypersurface in hyperbolic space is a geodesic sphere. Conversely, every geodesic ball attains equality. Inequality \eqref{eq:staticquadraticinequality} therefore holds under strict static convexity alone, while the convex normal-body construction in \autoref{lem:weakmixed} additionally assumes $h\geq0$.

For the local support variation $u_s=u+s\sqrt{1+u^2}$, the initial normal speed is $f=1$ and its derivative at fixed Euclidean Gauss direction is $-V_\nu/V$. In particular, the local variation satisfies
\begin{equation}\label{eq:staticBMdeficit} \left.\frac{d^2}{ds^2}I_{-1}(\Omega_s)^{1/(n+1)}\right|_{s=0}=-nI_{-1}^{-(2n+1)/(n+1)}\bigl(I_0^2-I_{-1}I_1\bigr). \end{equation}
Compare \eqref{eq:staticBMdeficit} with the parallel-domain identity \eqref{eq:parallelBMdeficit} below.

\subsection{Parallel domains}

For $h>0$, the outer parallel boundaries remain smooth, and \eqref{eq:HuLiweakparallel} shows that static convexity, including its strict form, is preserved.

\begin{proposition}\label{prop:intrinsicBM}
Let $\Omega$ be a compact domain with smooth connected boundary satisfying $h>0$ and $Vh-V_\nu g\geq0$, and let $\Omega_t$ be its outer parallel domain at distance $t\geq0$. For $0\leq j\leq n-1$,
\begin{equation}\label{eq:parallelvariations} I_{-1}^\prime=(n+1)I_0,\quad I_n^\prime=(n+1)I_{n-1},\quad I_j^\prime=(j+1)I_{j-1}+(n-j)I_{j+1}. \end{equation}
If $F(t)=I_{-1}(\Omega_t)^{1/(n+1)}$, then, with all integrals evaluated on $\Omega_t$,
\begin{equation}\label{eq:parallelBMdeficit} F''(t)-F(t)=-nI_{-1}^{-(2n+1)/(n+1)}\bigl(I_0^2-I_{-1}I_1\bigr)\leq0. \end{equation}
\end{proposition}

\begin{proof}
Use the weighted-volume first variation and \eqref{eq:intrinsicvariation} with $f=1$. By the divergence theorem, $\int_{\Sigma_t}V_\nu\dd\mu_t=I_{-1}(\Omega_t)$, while the hyperbolic Minkowski identities \cite[Eq.~(2.13)]{HuLi} imply $\int_{\Sigma_t}V_\nu E_j\dd\mu_t=I_{j-1}(\Omega_t)$ for $1\leq j\leq n$. This proves \eqref{eq:parallelvariations}. Differentiating $F$ and substituting $I_{-1}''=(n+1)(I_{-1}+nI_1)$ proves \eqref{eq:parallelBMdeficit}. The nonpositivity in \eqref{eq:parallelBMdeficit} follows from the weighted-volume inequality \eqref{eq:volume} on $\Omega_t$.
\end{proof}

\section*{AI disclosure}
AI tools were used during the development of this paper, in particular to assist with the geometric interpretation of the main theorem.

\end{document}